\documentclass[11pt]{article}

\newtheorem{theorem}{Theorem}[section]
\newtheorem{definition}[theorem]{Definition}
\newtheorem{prop}[theorem]{Proposition}
\newtheorem{example}[theorem]{Example}
\newtheorem{remark}[theorem]{Remark}

\newcommand{\simL}{\stackrel{\textit{l}}{\sim}}
 \newcommand{\TsfPair}{\underline{\underline{\mbox{Tsf Pair($K$)}}} }
\newcommand{\LinRel}{\underline{\underline{\mbox{Lin Rel($K$)}}} }
\newcommand{\simR}{\stackrel{\mathrm{r}}{\sim}}

\newcommand{\GA}{\Gamma}

\newcommand{\G}[2]{\ensuremath{\gamma_{#1,#2}}}  

\newcommand{\ddim}{\ensuremath{\mathrm{dim}}}

\begin{document}
\title{ Algorithms for computing linear invariants of directed graphs}
\author{Jacob Towber\\DePaul University\\Chicago,Illinois}
\date{\today}
\maketitle
\begin{abstract}
The present paper presents without proof, some algorithms which
solve the problem next to be explained. $($The proofs occur in a
longer version, submitted 3-11-2005 to \emph{Linear Algebra and
its Applications}(Elsevier)---the author will be happy to E-mail
this longer version in response to requests sent to
jtowber@uic.edu $)$
\smallskip

\noindent  Call two pairs $(M,N)$ and $(M',N')$ of $m\times n$
 matrices over a field $K$,  \emph{simultaneously
K-equivalent} if there exist square invertible matrices $S,T$ over
K, with $M'=SMT$ and $N'=SNT$. Kronecker \cite{Kronecker} has given a complete
set of invariants for simultaneous equivalence of pairs of
matrices.

Associate in the natural way to a finite directed graph $\Gamma$,
with $v$ vertices and $e$ edges, an ordered pair $(M,N)$ of
$e\times v$ matrices of zeros and ones. It is natural to try to
compute the Kronecker invariants of such a pair $(M,N)$,
particularly since they clearly furnish isomorphism-invariants of
$\Gamma$. Let us call two graphs `linearly equivalent' when their
two corresponding pairs are simultaneously equivalent.

     The purpose of the present paper, is to compute directly these
Kronecker invariants of finite directed graphs, from elementary
combinatorial properties of the graphs. A pleasant surprise is
that these new invariants are purely {\bf rational} ---indeed,
{\bf integral}, in the sense that the computation needed to decide
if two directed graphs are linearly equivalent
  only involves counting vertices in various finite graphs constructed
from each of the given graphs---
 and { \bf does not involve finding the irreducible
factorization of a polynomial over K} (in apparent contrast both
to the familiar invariant-computations of graphs furnished by the
eigenvalues of the connection matrix, and to the isomorphism
problem for general pairs of matrices.)
\end{abstract}

\section{Introduction: First Statement of the Problem}\label{S:intro}
    Let $\Gamma$ be a finite directed graph, with $v$ vertices and $e$
edges. Let us call $\Gamma$ {\bf reduced} if it contains no
`parallel' pairs of edges. (We allow `reduced' graphs to contain
`loops', i.e. $\Gamma$ may contain an edge which connects a vertex
to itself.) (The invariants to be defined, will turn out basically
to depend only on the `reduced' version of $\Gamma$. For the time
being, however, we do not assume that $\Gamma$ is reduced.)

If we choose (arbitrarily) orderings $<_E=(E_1,E_2,\cdots,E_e)$ of
the edges of $\Gamma$, and $<_V=(V_1,V_2,\cdots,V_v)$ of the
vertices, we then obtain two $e\times
v$ matrices $$M=M(\Gamma,<_E,<_V),N=N(\Gamma,<_E,<_V)$$ of 0's and
1's, given by:
\begin{equation}
M_{i,j}=\left\{{{1\mbox{ if $V_j$ is the initial vertex of $E_i$}}
\atop {0 \mbox{ otherwise}}} \right.
 \end{equation}
and similarly
\begin{equation}
 N_{i,j}= \left\{{{1\mbox{ if $V_j$ is the
terminal vertex of $E_i$}}
  \atop {0 \mbox{ otherwise}}}\right.
\end{equation}

Now suppose we are given a second directed graph $\Gamma^{'}$,
with the same numbers $v$ of vertices, and $e$ of edges as
$\Gamma$. Let $\Gamma^{'}$ be similarly associated with a pair
$(M', N')$ of $ e\times v$ matrices.

 Clearly,the two following assertions are equivalent:
\begin{flushleft}
Ia)\quad $\Gamma$ and $\Gamma'$ are isomorphic as directed graphs.
\smallskip

Ib) There exist $e\times e$ (resp.$v\times v$) permutation
matrices $S$ (resp. $T$) such that
\end {flushleft}
\begin{equation}
M'=SMT,\mbox{ and }N'=SNT.
\end{equation}
 The purpose of the present paper is to consider what happens
 when we replace the question of when
$\Gamma$ and $\Gamma'$ are related in the manner just described (a
question which seems immensely difficult),
  by the much easier modification of this question, suggested by
  the following two definitions.

  Let $K$ be any field.
\begin{definition}
Let $M,N,M',N'$ be matrices over the ground-field $K$, all of the
same size $m\times n$. We shall say that $(M,N)$ and $(M',N')$ are
\underline{simultaneously} \underline{ K-equivalent}, if there
exist invertible square matrices $S$ and $T$ over $K$, of sizes $m
\times m$ and $n\times n$ respectively, such that $$M'=SMT\,\mbox{
and } N'=SNT$$

\end{definition}
\begin{definition}\label{def:linequiv}
Let $\Gamma$ and $\Gamma'$ be finite directed graphs, with the
same number $v$ of vertices, and the same number $e$ of edges. Let
the pair $(M,N)$ of $e\times f$ matrices of 0's and 1's be related
to $\Gamma$ as described above (i.e. via 1),2) together with an
arbitrary choice of the orderings $<_E$ and $<_V$). Let $(M',N')$
be similarly associated with $\Gamma'$. Then we say that $\Gamma$
and $\Gamma'$ are \underline{$K$-linearly} \underline{equivalent},
if $(M,N)$ is simultaneously $K$-equivalent to $(M',N')$.
\end {definition}

We note that Def.\ref{def:linequiv}  is independent of the choice of vertex-
 and edge-orderings on $\Gamma$ and $\Gamma'$. We also note that
 the only difference between graph-isomorphism and $K$-linear
 equivalence, is that S,T in Eqn.3) are required to be
 permutation-matrices for $\Gamma$ and $\Gamma'$ to be isomorphic,
  while they only need to be invertible matrices over $K$ for
  $K$-linear equivalence to hold.

  Thus $K$-linear equivalence is a coarser equivalence relation than isomorphism;
  and the purpose of the present paper is to  establish some easily
  computed combinatorial invariants of directed graphs, which furnish
  what the author hopes is an
  efficient decision-procedure for $K$-linear equivalence.

\smallskip
The author would like to thank the University of Illinois at
  Chicago for its hospitality during
  the final work on this paper, and Schmuel Friedland and Raphael Loewy for
 encouraging conversations on these matters.
The author is honored to
  acknowledge his debt to Kronecker's fundamental 19th century
  work \cite{Kronecker}
  (extending earlier work of Weierstrass \cite{Weierstrass}) ,which is next to be
  discussed.
Thanks to my wife Diane, for her good humour which made
  these efforts practicable...

  \section{A More Detailed Statement of the Problem}\label{S2}
Let K be a given field.

 We next wish to review the details of
Kronecker's canonical form, for pairs of matrices over $K$ under
simultaneous $K$-equivalence (cf.Def.1.2). The computational
details involved in reduction to Kronecker's canonical form, seem
to the author to become clearer,\  if stated in terms of
vector-spaces over $K$ and $K$-linear transformations, rather than
in terms of matrices. Thus, let us consider the problem of
classifying pairs of maps from one vector-space over $K$ into a
second--- i.e. of classifying diagrams
\begin{equation}\label{E:diagram}
E{\stackrel{\mu}{\longrightarrow} \atop
\stackrel{\nu}{\longrightarrow} }  V
\end{equation}
where $E$ and $V$ are finite-dimensional vector-spaces over the
ground-field $K$, and where $\mu$ and $\nu$ are $K$-linear
transformations---such a diagram will be referred to as a
{\em transformation-pair} over $K$.
 We shall also use the notation
$$[\mu,\nu:E\rightarrow V]$$ to refer to the transformation-pair (\ref{E:diagram}).

Of course, to have a well-specified classification problem, we
must specify precisely what meaning we wish to attach to `isomorphism'
 between two such diagrams
 (\ref{E:diagram}).  We do this in the obvious way: by a {\bf
$K$-isomorphism} from $[\mu,\nu:E\rightarrow V]$ to
$[\mu',\nu':E'\rightarrow V']$ will be meant an ordered pair
$$(\alpha:E\to E',\beta:V\to V')$$ of $K$-linear isomorphisms,
such that
$$\mu'\circ\alpha=\beta\circ\mu,\,\mbox{ and
}\nu'\circ\alpha=\beta\circ\nu$$
(This is just another way of describing the earlier notion of
simultaneous equivalence of pairs of matrices, the matrices being those that
represent $\mu,\nu$ with respect to a choice of $K$-bases for $E$ and $V$.)
\smallskip
\par Thus one is led to study the classification problem for the
category whose objects are diagrams (4), and whose morphisms are
the $K$-isomorphisms just  specified---let us denote this category
by \underline{\underline{Tsf Pair($K$)}} \,. One slight further
delay, before we finally get to Kronecker's beautiful
classification of diagrams (\ref{E:diagram}): let us  note how all
this is connected ( by a process of `linearization' ) with
directed graphs:

\begin{definition}\label{def:KLin}
Let $\Gamma$ be a finite directed graph, with vertex set $V$ and
whose set of edges is $E$. We then construct a corresponding
element
$$K(\Gamma)=[\mu_{\Gamma} ,\nu_{\Gamma }:E_{\Gamma }\rightarrow
V_{\Gamma }]$$ in \TsfPair as follows:

We take $E_{\Gamma }$ to be the vector-space over $K$ freely
generated by the basis $E$, $V_{\Gamma } $ to be similarly $K$-free on the
set V of vertices, while the $K$-linear transformations
$$\mu_{\Gamma},\nu_{\Gamma}:E_{\Gamma }\rightarrow V_{\Gamma }$$
are defined on the basis-vectors $e\in E$ for $E_{\Gamma }$, by:
$$\mu_{\Gamma}(e)=\mbox{initial vertex of e},\;\,\nu_{\Gamma}(e)=\mbox{terminal vertex of e.}$$
We shall refer to this transformation-pair $K(\Gamma)$ in
\TsfPair as the {\bf $K$-linearisation}
of $\Gamma$.
\end{definition}
\smallskip
\noindent \underline{{\bf EXAMPLE 2.1:}} Consider the directed
graph $\Gamma_1$ furnished by the following diagram:
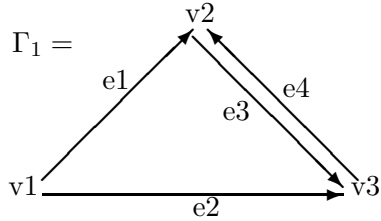
\begin {center}
 \setlength{\unitlength}{1cm}
 \begin{picture}(4,2.5)
 \thicklines
\put(0,0){\vector(1,0){4}} \put(0,.2){\vector(1,1){2}}
\put(2,2.1){\vector(1,-1){2}} \put(4.2,.2){\vector(-1,1){2}}
\put(-.45,0){v1} \put(1.9,2.3){v2} \put(4.1,0){v3} \put(.8,1.4){e1}
\put(2,-.28){e2}
 \put(2.4,1){e3} \put(3.2,1.3){e4}
 \put(- .4,1.9){$\GA_1=$}
\end{picture}
\end{center}
\medskip
Then the $K$-linearisation of $\Gamma_1$ is [$\mu,\nu:E\to V$],
 where
 $$E=K\{e_1,e_2,e_3,e_4\},\;V=K\{v_1,v_2,v_3\}$$
 and where (all x's lying in $K$)
\begin{eqnarray}\label{E:example1}
 \mu(x_1e_1+x_2e_2+x_3e_3+x_4e_4)&=&(x_1+x_2)v_1+x_3v_2+x_4v_3\\
\nu (x_1e_1+x_2e_2+x_3e_3+x_4e_4)&=&(x_1+x_4)v_2+(x_2+x_3)v_3\nonumber
\end{eqnarray}
 \medskip
 \par We now proceed to sketch Kronecker's classification results,
 stated in terms of the category \TsfPair.
 These results will here be stated without proof (for which see
\cite{Gantmacher},Chap.XII;
 \cite{Turnbull-Aitken},Chap.IX; or, \cite{Kronecker}. Also, in \cite{Towber}
 a direct proof is given for the classification of the category \LinRel
 defined in \S\ref{S:7} below; it is then easy (using the
 results obtained in \S\ref{S:7}) to deduce the classification of
 \TsfPair from that of \LinRel.)

 \bigskip
 \par In the first place, there is an obvious notion of direct sum
 defined as follows on the category
 \TsfPair\,: the direct sum
 $$[\mu,\nu:E\to V]\oplus [\mu',\nu':E'\to V'] $$
of two objects, is the object
$$[\mu\oplus\mu',\nu\oplus\nu':E\oplus E'\to V\oplus V']$$
There is a unique {\bf\it zero object} in \TsfPair, namely
$$[0,0:0\to 0]$$

 An object $[\mu,\nu:E\to V]$ in \TsfPair will be called {\bf\it
 indecomposable} if it is not the zero object, and is not
 isomorphic to the direct sum of two non-zero
 objects in \TsfPair.
\par Kronecker's theory then furnishes us with the following five
families of indecomposable objects in
\TsfPair, such that every indecomposable
object in this category, is isomorphic to a unique object in this
list:
 \smallskip
\par In this listing, it will be convenient to denote by
 $\{f_1^{(n)},\cdots,f_n^{(n)}\}$ the usual standard basis for $K^n$, so that
 $$f_1^{(n)}=\overbrace{(1,0,\cdots,0)}^{\mbox{$n$ entries}},f_2^{(n)}=(0,1,\cdots,0),\cdots,f_n^{(n)}
 =(0,0,\cdots,1)$$\par
 Also, we shall denote by $I_n$ the identity map on $K^n$.\par
 \smallskip
\noindent\underline{{\bf TYPE $^0T_n$}}\ \ $(n=1,2,3,\cdots)$

 This is the object $[\mu_n,I_n:K^n\to K^n]$ where ${\mu}_n$ is the
 nilpotent $K$-linear transformation on $K^n$ which maps:
 \begin{equation}\label{eq:nilp}
 \mu_n: f^{(n)}_1\mapsto f^{(n)}_2\mapsto\cdots\mapsto f^{(n)}_n\mapsto 0
 \label{eq:nilp}
 \end{equation}
\underline{{\bf TYPE $T_n^0$}}\ \ $(n=1,2,3,\cdots)$ \par This is
the object $[I_n,\mu_n:K^n\to K^n]$, where $\mu_n$ is
 still given by (\ref{eq:nilp}).
\underline{{\bf TYPE $T_n$}}\ \ $(n=0,1,2,\cdots)$\par This is the
object $[\kappa_n,\lambda_n:K^n\to K^{n+1}]$, where
  $\kappa_n,\lambda_n$ are defined by:
$$\kappa_n(f^{(n)}_1)=f^{(n+1)}_1,\kappa_n(f^{(n)}_2)=f^{(n+1)}_2,\cdots,
\kappa_n(f^{(n)}_n)=f^{(n+1)}_n$$ and
$$\lambda_n(f^{(n)}_1)=f^{(n+1)}_2,\lambda_n(f^{(n)}_2)=f^{(n+1)}_3,\cdots,
\lambda_n(f^{(n)}_n)=f^{(n+1)}_{n+1}$$

 (When $n=0$, this is to be understood as the object $[0,0:0\to
  K]$)
 \underline{{\bf TYPE $^0T_n^0$}}\ \ $(n=0,1,2,\cdots)$\par This is
the object $[\kappa_n',\lambda_n':K^{n+1}\to K^n]$, where
  $\kappa_n',\lambda_n'$ are defined by:
$$\kappa_n'(f^{(n+1)}_1)=f^{(n)}_1,\kappa_n'(f^{(n+1)}_2)=f^{(n)}_2,\cdots,
\kappa_n'(f^{(n+1)}_n)=f^{(n)}_n\mbox{, but }
\kappa_n'(f^{(n+1)}_{n+1})=0\;;$$ and
$$\lambda_n'(f^{(n+1)}_1)=0,\lambda_n'(f^{(n+1)}_2)=f^{(n)}_1,
\lambda_n'(f^{(n+1)}_3)=f^{(n)}_2,\cdots,\lambda_n'(f^{(n+1)}_{n+1})=f^{(n)}_n$$

 (When $n=0$, this is to be understood as the special object
 $[0,0:K\to 0]$ ---this special object plays a somewhat exceptional role, at
 several places in the investigations that follow.)

\noindent \underline{{\bf TYPE $S(p(X)^n$}})
 Here $n$ denotes a positive integer, and $p(X)$---subject to the
 condition $p(X)\neq X$---denotes a monic irreducible polynomial
 in the polynomial ring $K[X]$ in one indeterminate $X$ over $K$.
 $S(p(X)^n)$ is then defined to be the object $[\xi,I(V):V \to V]$
where $V=V(n,p)=K[X]/(p(X))^n$, where  $I(V)$ is the identity map on
V, and where $\xi$ is the $K$-endomorphism of $V$ given by
multiplication by $X$.\par \medskip Let us denote by
$\mathcal{L}(K)$ the list (just described) of indecomposable
objects in \TsfPair .

\bigskip
\begin{remark}\label{SF}
\end{remark}

\noindent More generally, for any monic
polynomial $F(X)\in K[X]$, it will be quite convenient for us to denote by
$S(F(X))$ the transformation-pair
$$[\phi,I:K[X]/(F)\to K[X]/(F)]\;,$$
 where $\phi$
is induced via multiplication by $X$ (so that $F(\phi)=0$). Note that, in particular,
$S(X^n)$ is $^0T_n$.

Of course, $S(F(X))$ is not in general indecomposable in \TsfPair;
rather, if F has the irreducible factorization in $K[X]$ given by
\begin{equation}\label{E:CRT}
F=\prod_{i=1}^n p_i(X)^{n_i}\mbox{ ,then  }S(F(X))=
\bigoplus_{i=0}^n S(p_i(X)^{n_i})
\end{equation}
(In our decomposition in \TsfPair of the $K$-linearizations of graphs, we shall
begin by splitting off certain of these $S(F)$, and then make use of (\ref{E:CRT}))

 In our present formulation, Kronecker's theory asserts that
every object $\pi$
 in \TsfPair is isomorphic to a direct sum of
 a finite number of indecomposable objects; and asserts that, given an
 indecomposable object $\kappa$ in the above list
 $\mathcal{L}(K)$,
 the multiplicity $[\pi:\kappa]$
 with which $\kappa$ occurs in such a direct sum decomposition for
 $\pi$, is the same for all such decompositions. Thus, a complete
 list of invariants for objects in
 \TsfPair is afforded by the collection of
 these multiplicities
 \begin{equation}\label{E:inv}
\{[\pi:\kappa] |\kappa\in \mathcal{L}(K)\}
 \end{equation}
 (which are non-negative
integers; for given $\pi$, all but finitely many of these
integers vanish.)\par
\medskip
We shall also use the following notation for the invariants
associated to the first four infinite families in the above list $
\mathcal{L}(K)$: for $\pi$ any transformation-pair over $K$, we
shall set \par
\begin{equation}\label{E;tDef}
\left \{
\begin {array} {c}
t_n(\pi):=[\pi:T_n]=\mbox{ multiplicity of $T_n$ in }\pi\quad (0\leq n)\\
\mbox{ and similarly} \\
^0t_n(\pi):=[\pi:^0\!\!T_n]\quad (1\leq n) \\
t_n^0(\pi):=[\pi:T_n^0] \quad (0\leq n) \\
^0t_n^0(\pi):=[\pi:^0\!\!T_n^0] \quad (0\leq n)
\end{array}
\right.
\end{equation}
\medskip
\noindent \underline{{\bf A FURTHER HISTORICAL NOTE:}}

\noindent The work of Kronecker just sketched, was based on the
earlier classification results in Weierstrass' 1867 paper
\cite{Weierstrass}.
 Restating
Weierstrass' work in terms of transformation-pairs: in it were
classified
 those transformation-pairs $$[\mu,\nu:E\to V]$$ which are
 \underline{{\bf regular}}, i.e. have both the following
 properties:
 $$ \mbox{1) }\ddim(E)=\ddim(V), \mbox{ and 2) } \exists c\in K
 \mbox{ so } |\mu+c\nu|\neq 0$$

Namely, Weierstrass showed (in our present terms) that every
regular transformation-pair, is uniquely isomorphic to a direct
sum of the indecomposable types
$$^0T_n \mbox{ and }T^0_n\quad(n>1);\qquad S(p^n)\qquad
 p\;\mbox{irreducible,}n>0$$
defined above.

 This work of Weierstrass is explained in Gantmacher(\cite{Gantmacher}), pp.24-28, and also
 in Turnbull-Aitken(\cite{Turnbull-Aitken}), pp.113--118
 ---both these  books go from there to explain the extension of
 this work, in Kronecker's 1890 paper \cite{Kronecker}, to the full category
 $\TsfPair$ (as explained above) which involved adding the two
 infinite classes of non-regular indecomposable forms
 $$T_n \mbox{ and }^0T^0_n\qquad(n\geq 0)\;.$$

\medskip
\noindent \underline{{\bf EXAMPLE 2.1 RE-VISITED:}}\par

 In Example 2.1 we computed the $K$-linearization of the directed graph
$\Gamma_1$ pictured above. Let us, in this simple example,
compute the associated Kronecker invariants:\par
               We begin by looking inside $K(\Gamma_1)$ for
a copy of Type $^0T_1^0$, i.e. of the indecomposable
transformation-pair $[\mu,\nu :K^2\to K]$ indicated by the
following picture:
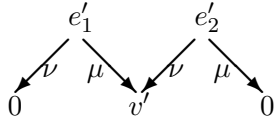
\begin {center}
 \setlength{\unitlength}{1cm}
 \begin{picture}(4,2.5)
 \thicklines
\put(0,2){\vector(-1,-1){.7}}
\put(.2,2){\vector(1,-1){.7}}
 \put(1.68,2){\vector(-1,-1){.7}}
 \put(1.88,2){\vector(1,-1){.7}}
\put(0,2.2){$e'_1$} \put(1.68,2.2){$e'_2$} \put(-.8,1 ){$0$}
\put(.8,1){$v'$}  \put(2.55,1){$0$} \put(-.35,1.5){$\nu$}
\put(.25,1.5){$\mu$} \put(1.33,1.5){$\nu$} \put(1.93,1.5){$\mu$}
\end{picture}
\end{center}
Using  equations (\ref{E:example1}), we see first that $\mu
e'_2=0$ implies that $e'_2$ is a scalar multiple of $e_1-e_2$, say
$$e'_2=e_1-e_2,v'=\nu e_2'=\nu (e_1-e_2)=v_2-v_3;.$$
Then $$\mu e_1'=v'=v_2-v_3,\nu e_2'=0$$ specify $e_1'$ uniquely as
$e_1-e_2+e_3-e_4$ This shows the Kronecker invariant
$^0t^0_1(\Gamma_1)=[K(\Gamma_1):\ ^0T^0_1]$ has the value $1$.

We next look inside $K(\Gamma_1)$ for a copy of Type $T^0_1$. Thus, we try to find
$$e''\in E=K\{e_1,e_2,e_3,e_4\}, v''\in V=\{v_1,v_2,v_3\}$$
such that $\mu e''=v''$ and $\nu e''=0$. Using equations (\ref{E:example1}),
it is readily seen that the most general solution to these equations
is given by
$$e''=\alpha(e_1-e_4)+\beta(e_2-e_3),
v''=\alpha(v_1-v_3)+\beta(v_1-v_2)$$ (Almost any choice of
$(\alpha, \beta)$ is satisfactory; choices which will fail to give
a Kronecker decomposition are those for which $\beta=-\alpha$, so
that $e''$ is a scalar multiple of $e_2'$ and our second object is
contained inside the first.) Let us choose $\alpha=1,\beta=0$, so
our copy of Type $T^0_1$ has
$$e''=e_1-e_4,\; v''=v_1-v_3$$

Finally, it is readily verified that a copy of S(X-1) inside $K(\Gamma_1)$, i.e.
a solution to $\mu e'''=\nu e'''=v'''$, is furnished by
$$e'''=e_1-e_2+e_3,\; v'''=v_2$$
It is easy to check that (with the choices indicated above)
$\{e'_1,e'_2,e'',e'''\}$ are linearly independent over $K$, as are
$\{v',v'',v'''\} $, so that $K(\Gamma_1)$ is the direct sum of the
three sub-objects constructed (which could however have been
selected for this purpose in infinitely many other ways), and we
have:
\begin{equation}\label{decomp:G1}
 K(\Gamma_1)\simeq\ ^0T^0_1 \oplus T^0_1\oplus S(X-1)
\end{equation}
{\bf EXAMPLE 2.2:} Let us consider a second example (whose
properties will be useful in the next section). Namely, let
$S_2=\{v_1,v_2\}$ be a set with two elements, and let $\GA_2$ be
the directed graph, on these two vertices, which contains
precisely one edge $e_{i,j}$ going from $v_i$ to $v_j$ for
$$i,j\in \{1,2\}$$
The $K$-linearization is then
$$K(\GA_2)=[\mu_2,\nu_2:E_2\to V_2]$$
where $V_2$ is $K$-free on $v_1,v_2$, $E_2$ is $K$-free on the four
edges
$$e_{1,1},\,e_{1,2},\,e_{2,1},\,e_{2,2}$$
and where
\begin{eqnarray}
\mu_2(a_{1,1}e_{1,1}+a_{1,2}e_{1,2}+a_{2,1}e_{2,1}+a_{2,2}e_{2,2})
&=&(a_{1,1}+a_{1,2})v_1+(a_{2,1}+a_{2,2})v_2  \nonumber \\
\nu_2(a_{1,1}e_{1,1}+a_{1,2}e_{1,2}+a_{2,1}e_{2,1}+a_{2,2}e_{2,2})
&=&(a_{1,1}+a_{2,1})v_1+(a_{1,2}+a_{2,2})v_2   \nonumber
\end{eqnarray}
If
$$\alpha=a_{1,1}e_{1,1}+a_{1,2}e_{1,2}+a_{2,1}e_{2,1}+a_{2,2}e_{2,2}$$
satisfies $\mu_2(\alpha)=\nu_2(\alpha)=0$,then it follows that
$$a_{2,1}=-a_{1,1},a_{1,2}=-a_{1,1},a_{2,2}=a_{1,1}\;,$$
so $Ker\; \mu_2\cap Ker\; \nu_2$ consists of multiples of
$$e_{1,1}-e_{1,2}-e_{2,1}+e_{2,2}\;.$$
This gives a copy inside $K(\GA_2)$ of
$$^0T^0_0=[0,0:K\to 0]$$
unique up to scalar multiple.The general results to be obtained in
the sequal, will imply the Kronecker decomposition
\begin{equation}\label{decomp:G2}
K(\GA_2)\cong\, ^0\!T^0_0\oplus ^0\!T^0_1\oplus S(X-1)
\end{equation}
The reader may wish to prove this, directly from the definitions.
\medskip
\par We now conclude this section by re-stating (in sharper focus)
the question, sketched in the preceding section, (and discussed
for the two preceding special cases),
 whose general solution will occupy the remainder of this
paper:

\par {\bf Given a directed graph $\Gamma$, how
can we efficiently compute the Kronecker invariants}
 $$\{\;[K(\Gamma):\kappa]\;|\;\kappa \in \mathcal{L}(K)\}$$
 {\bf of the $K$-linearization of $\Gamma$ }?

 \section{Reduction to the Case of Binary Relations}\label{S3}
 Our next step is an easy preliminary reduction of the problem just
 stated:

 Namely, let $\Gamma$ be a finite directed graph; denote by
 $E_{\Gamma}$ resp. $V_{\Gamma}$ the set of edges resp. vertices of
$\Gamma$.
Let $\sim_{\Gamma}$ denote the equivalence-relation of
`parallelism' on the set of edges of $\Gamma$, so that if $e$ and
$e'\in E_{\Gamma}$:
\begin{flushleft}
$e \sim_{\Gamma} e' \Longleftrightarrow e$ and $e'$ have the same
initial vertex, and also the same final vertex.
\end{flushleft}
\smallskip
\par In \S1 we  defined $\Gamma$ to be \emph{ reduced} if $\Gamma$ does not
contain a pair $e,e'$ of distinct edges which are parallel. Even
if $\Gamma$ does contain parallel edges, the following obvious
construction produces a finite directed graph $\Gamma^{red}$ which
is reduced:

We define the set
$$E_{\Gamma^{red}}:=E_{\Gamma}/\sim_{\Gamma}$$
of edges of $\Gamma^{red}$ to be the set of
$\sim_{\Gamma}$-equivalence-classes; we also set
 $$V_{\Gamma^{red}}:=V_{\Gamma}$$
 i.e., $\Gamma^{red}$ is to have the same vertices as $\Gamma$.
 Finally, if $e$ is an edge of $\Gamma$, and if cls $e \in
 E_{\Gamma^{red}}$ is the $\sim_{\Gamma}$-equivalence-class
 containing $e$, then the initial (resp. terminal) vertex in
 $ \Gamma^{red}$ of cls $e$ are defined to coincide with the
 initial (resp. terminal) vertex in $\Gamma$ of $e$.

Let us call $\Gamma^{red}$ the \emph{reduced form} of $\Gamma$.
With the notation just explained, we have:

\begin{prop}\label{prop:red}
Let the equivalence-relation $\sim_{\Gamma}$ partition
$E_{\Gamma}$ into $s$ equivalence classes, of cardinalities
 $n_1,n_2,\cdots, n_s$. Then the $K$-linearization
 $K(\Gamma)$ of $\Gamma$, is the direct sum of $K(\Gamma^{red})$ with
 \begin{equation}\label{E:mult}
 n_1+n_2+\cdots +n_s-s
 \end{equation}
 copies of the irreducible type $^0T_0^0=[0,0:K\to 0]$.
 \end{prop}
\par\noindent {\bf PROOF:} Let $E_1,\cdots,E_s$ be the
$\sim_{\Gamma}$-equivalence classes, and let
$$E_i=\{e^{i}_1,\cdots,e^i_{n_i}\}\mbox{ for }1\leq i\leq s\;.$$
Then set
$$\overline{E}:=\{e^1_1,e^2_1,\cdots,e^s_1\}$$
and
$$f^i_j:=e^i_j-e^i_1\in E_{\Gamma}\mbox{ for }1\leq i\leq s,2\leq
j\leq n_i\;.$$ It is readily verified that
$\overline{\Gamma}:=(\overline{E},V)$ forms a sub-graph of
$\Gamma$ isomorphic to $\Gamma^{red}$, and that $K(\Gamma)$ is the
direct sum of $K(\overline{\Gamma})$ with the $\sum (n_i-1)$
sub-objects
$$[0,0:Kf^i_j\to 0]\cong ^0\!\!T_0^0  \qquad (1\leq i\leq s,2\leq
j\leq n_i)\;.$$ This proves the asserted proposition.\par
\bigskip
\noindent {\bf CAUTION:} It does \emph{not} follow from the Prop.
\ref{prop:red} just proved, that the multiplicity
$^0t^0_0\;(K(\GA))$ with which $^0T^0_0$ occurs as a direct
summand in $K(\GA)$ is given by (\ref{E:mult})---because also
$K({\GA}^{red})$ may contain one or more summands of type
$^0T^0_0$. (See Example 2.2 above for an example of a reduced
graph whose $K$-linearization contains
 a direct summand of type $^0T_0^0$.)

\medskip
 Because of  the proposition Prop. \ref{prop:red} just proved, we may as
well restrict to reduced directed graphs. If $\Gamma$ is a reduced
directed graph, with vertex-set $V$, we may as well consider the
set $E$ of edges to be a sub-set of $V\times V$,\\\underline{ and this is the
viewpoint which we shall adopt for the remainder of this}\\\underline{ paper.}
Then $E\subseteq V\times V$ endows our set-up with the familiar structure
of a binary relation $E$ on a finite set $V$.

This modified point of view, is a bit better adapted to the
combinatorial constructions in the next section. It makes
available the various basic operations on the set of binary
relations on a given set $S$---e.g., given binary operations $R,
R'$ on $S$, we may form the following further binary relations on
$S$:
$$R\cup R',R\cap R',\mbox{ the converse relation }R^{-1},
 \mbox{ and the composite }RR'\;.$$
---operations which will make possible the constructions of the next section.

\smallskip
 Just to clarify our terminology: for the rest of
 this paper, a {\bf binary relation} is an ordered pair
$$\GA=(R,V)\mbox{ where }R\subseteq V\times V\;.$$
---{\bf in the present work, ``binary relation'' is never taken in
the more general sense of a subset $R\subseteq S\times S'$ with $S$ and
$S'$ distinct sets.}

$\GA$ is \emph{finite} if $V$ is a finite set. We shall sometimes vary
our language by writing instead
$$\GA=(R\subseteq V\times V)$$
 For $v$ and $v'$ in $V$, the notations $(v,v')\in R$ and $vRv'$
 are synonymous.

\bigskip
Let us spell out how the problem being studied looks with this minor change in
language:

Let $\GA=(R,V)$ denote a binary relation on a finite
set $V$. This corresponds to the reduced directed graph, whose
edge-set coincides with $R$, whose vertex-set coincides with $V$,
and such that the ordered pair $(v,v')$ in $R$ is regarded as an
edge, with initial vertex $v$, and terminal vertex $v'$. The
 $K$-linearization $K(\GA)$ of this directed graph consists of the
 transformation-pair
 $$[\mu_\GA,\nu_\GA:KR\to KV]\in \mbox{\TsfPair}$$
where $KV$ (resp. $KR$) is the vector-space over K free on the
basis $V$ (resp. $R$), and where the $K$-linear maps $\mu_\GA,\nu_\GA$
are defined on basis-vectors $(v,v')\in R$ by
$$\mu_\GA(v,v')=v,\;\nu_\GA(v,v')=v'$$
For each $\kappa\in \mathcal{L}(K)$ we shall denote by
$[\GA:\kappa]$ the multiplicity with which the indecomposable object
 $\kappa$ appears as a direct summand in $K(\GA)$.

The fundamental problem studied in this paper, may now be re-stated
as the computation, for finite binary relations, of these integers
$$\{[\GA:\kappa] | \kappa\in \mathcal{L}(K)\}$$
In the next section, we define some combinatorial constructions on
binary relations $R$, which will enable us to accomplish this.

\section{Left and Right Contractions of Binary Relations}
\label{S:4} Let $\GA=(R,S)$ denote a binary relation on a finite
set $S$. (We consider $S$ to be part of the given structure $\GA$,
i.e. if $S$ is a proper subset of $S'$, $(R,S)$ and $(R,S')$ count
as distinct binary relations.)
\begin{definition}\label{def:quotient}
Let $\sim$ be an equivalence relation on $S$.We then denote by
$S/\sim$ the set of $\sim$-equivalence classes, and by
$$\GA/\sim=(R/\sim, S/\sim)$$
the binary relation on $S/\sim$ specified by
$$\alpha(R/\sim)\beta\Longleftrightarrow \exists a\in
\alpha,b\in\beta \mbox{ with } aRb\mbox{ (for   }\alpha,\beta\in
S/\sim)$$ Also,we denote by $\pi(\GA,\sim)$ the canonical
epimorphism
$$\pi(\GA,\sim):S\to S/\sim$$ which maps each $s\in S$ into the
$\sim$-equivalence class containing s.
\end{definition}
\begin{definition}\label{def:contract}
We denote by
$$\stackrel{\textit{l}}{\sim}=\stackrel{\textit{l}}{\sim}(\GA) $$
the equivalence relation on $S$ generated by $R^{-1}R$, i.e. the
smallest equivalence relation on $S$ such that
$$xRy \mbox{ and }xRy' \Rightarrow y
\stackrel{\textit{l}}{\sim}y'\;;$$ and we define the {\bf left
contraction} of R to be the binary relation
 $$C_L(\GA):=\GA/\stackrel{\textit{l}}{\sim}\mbox{ on the set }
 V/\stackrel{\mathrm{l}}{\sim}$$ (as given by Def.\ref{def:quotient}).

Similarly, we denote by
$$\stackrel{\mathrm{r}}{\sim}=\stackrel{\mathrm{r}}{\sim}(\GA) $$
the equivalence relation on $S$ generated by $RR^{-1}$, i.e. the
smallest equivalence relation on $S$ such that
$$xRy \mbox{ and }x'Ry \Rightarrow x
\stackrel{\mathrm{r}}{\sim}x'\;;$$ and we define the {\bf right
contraction} of $\GA$ to be the binary relation
 $$C_R(\GA):=\GA/\stackrel{\mathrm{r}}{\sim}\mbox{ on the set }
 V/\stackrel{\mathrm{r}}{\sim}\;\,.$$
\end{definition}
\begin{remark}\label{remark:converse}
\end{remark}

If we denote the  {\bf converse relation} $(R^{-1},S)$ to $\GA$ by
$\GA^{-1}$, where (for all $s$ and $s'$ in $S$),
$$sR^{-1}s'\Longleftrightarrow s'Rs\;,$$
 then clearly
$$C_l(\GA^{-1})=(C_r\GA)^{-1},\mbox{ and }C_r(\GA^{-1})
=(C_l\GA)^{-1}\;.$$

 \bigskip
 \begin{example}\label{Ex4.4}
\end{example} Consider the  binary relation $\GA_3$ given by this diagram:
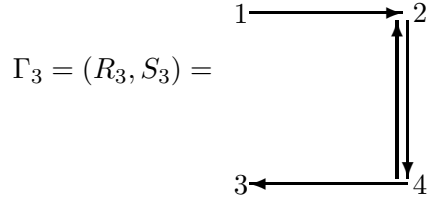
\begin{center}
 \setlength{\unitlength}{.7cm}
\begin{picture}(0,4)
\thicklines \put(3,0){\vector(-1,0){3}} \put(-.3,-.2){3}
\put(3.1,-.2){4} \put(2.8,0.1){\vector(0,1){3}}
\put(3,3.1){\vector(0,-1){3}}
 \put(0,3.25){\vector(1,0){2.9}}
\put(-.3,3.05){1} \put(3.1,3.07){2} \put(-4.5,2){$\GA_3=(R_3,S_3)=$}
\end{picture}
\end{center}
Then $4R_32,4R_33$ imply $2\simL 3$, and the $\simL$-equivalence
classes are $\{1\},\{2,3\},\{4\}$.

Thus the left contraction of $(R_3,S_3)$ has three vertices, and
may be drawn as
 \setlength{\unitlength}{.7cm}
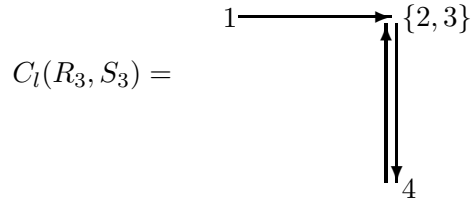
\begin{center}
\begin{picture}(0,4)
\thicklines \put(3.1,-.2){4} \put(2.8,0.1){\vector(0,1){3}}
\put(3,3.1){\vector(0,-1){3}}
 \put(0,3.25){\vector(1,0){2.9}}
\put(-.3,3.05){1} \put(3.1,3.07){$\{2,3\}$}
\put(-4.3,2){$C_l(R_3,S_3)=$}
\end{picture}
\end{center}

Similarly, $1R_32$ and $4R_32$ imply $1\simR 4$, and the right
contraction of $(R_3,S_3)$ is given by
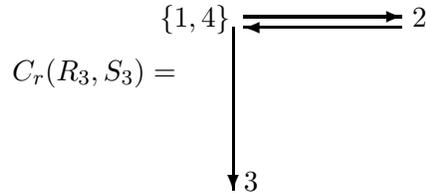
\begin{center}
 \setlength{\unitlength}{.7cm}
\begin{picture}(0,4)
\thicklines \put(-.1,3){\vector(0,-1){3.1}}
\put(-1.5,3){$\{1,4\}$} \put(.1,3.2){\vector(1,0){3}}
\put(3.1,3){\vector(-1,0){3}} \put(.1,-.1){3} \put(3.3,3){2}

 \put(-4.3,2){$C_r(R_3,S_3)=$}
\end{picture}
\end{center}

\section{Statement of the Main Theorem}\label{S:main}
The present section contains a number of statements, which
together furnish an algorithm which constitutes a complete
solution to the computational problem stated in Sections 1 and 2.
{\bf Note: No proofs are presented in this section, or the next,
in order not to interrupt the exposition of these
results}---results which will be stated in this section, and will
be illustrated by some specific computations in the next section.
 Following these two sections, the remainder of the present paper then contains the
proofs of these results.

In this section, $\GA=(R,S)$ will denote a finite binary relation
(i.e., a reduced directed graph.)
\begin{prop}\label{ClCrCommute}
Let $\GA=(R,S)$ be a binary relation; then there is a natural
graph-isomorphism\footnote{Actually, this may be strengthened, to
assert a closer identity between $C_lC_r\GA$ and $C_rC_l\GA$;
cf.Th.7.12. This stronger form is useful in checking computations
such as those in \S\ref{S:EX}.} between $C_lC_r\GA$ and
$C_rC_l\GA.$
\end{prop}
\begin{definition}\label{gamma}
For m,n any non-negative integers, we may define (using
the preceding Prop.(\ref{ClCrCommute}))
a binary relation
$$C_l^mC_r^n\GA=(^m\!R^n,^m\!\!S^n)$$
(which we take to be $\GA$ if $m=n=0$.) Then by
$\gamma_{m,n}(\GA)$ will be meant the cardinality of the set
$^m\!S^n$ .
\end{definition}

\begin{example}\end{example}
Let us re-examine Example 4.4 in Section \ref{S:4}. For the binary
relation $\GA_3$ in this example, we have
$$\G{0}{0}(\GA_3)=4,\G{1}{0}(\GA_3)=3,\G{0}{1}(\GA_3)=3\;.$$
It is readily verified that
$$C_l^n \GA_3\cong C_l \GA_3,\mbox{ and } C_l^n\GA_3\cong C_l\GA_3
\mbox{ for all } n\geq 1\;.  $$ whence
$\G{n}{0}(\GA_3)=\G{0}{n}(\GA_3)=3$ for $n\geq 1$.

 $C_lC_r\GA_3\cong C_rC_l$ is given by
 $$\{1,4\}{\longleftarrow \atop \longrightarrow} \{2,3\}$$
and since this graph is stable under the actions of $C_l$
and $C_r$, we see that
$$\G{m}{n}(\GA_3)=2 \mbox{ for }m\geq 1,n\geq 1\;.$$
\begin{theorem}\label{th:main}
$({\bf Main\; Theorem, Part\; One})$ Let $\GA=(R,S)$ be a finite
binary relation, with
$K$-linearization $K(\GA)$. Then:\\
for $n\geq 1$,
\begin{eqnarray}
^0t_n(K(\GA))&=&\G{n+1}{n-1}-\G{n}{n}-\G{n+1}{n}+\G{n}{n+1}  \label{E:ZT}\\
            &=&\G{n-1}{n+1}-\G{n}{n}-\G{n-1}{n}+\G{n}{n-1} \nonumber \\
t^0_n(K(\GA))&=&\G{n-1}{n+1}-\G{n}{n}-\G{n}{n+1}+\G{n+1}{n}  \label{E:TZ} \\
           &=&\G{n+1}{n-1}-\G{n}{n}-\G{n}{n-1}+\G{n-1}{n} \nonumber
\end{eqnarray}
For $^0t^0_n$---still with the restriction $n>0$---
 we must distinguish two cases, according as $n$ is even or odd:
\begin{eqnarray}
^0t^0_{2p+1}(K(\GA))&=&\G{p}{p}-\G{p}{p+1}-\G{p+1}{p}+\G{p+1}{p+1}\label{E:ZTZodd} \\
^0t^0_{2p}(K(\GA))&=&\G{p}{p-1}-\G{p}{p}-\G{p+1}{p-1}+\G{p+1}{p}\label{E:ZTZeven} \\
                 &=&\G{p-1}{p}-\G{p}{p}-\G{p-1}{p+1}+\G{p}{p+1}\nonumber
\end{eqnarray}
\end{theorem}

\bigskip
This still leaves open the multiplicities in $K(\GA)$ of the
indecomposable types $T_n$ and $S(n,p(X))$ --- also,
$^0t^0_0(\GA)$ remains open. To determine these, we first need the
following definition (and also, of course, proof of the facts
assumed by this definition---as already promised, all such proofs
will be supplied in subsequent sections of this paper.)
\begin{definition}\label{stable}
Let $\GA=(R,S)$ be a finite binary relation. Then
the graphs $C_l^mC_R^n \GA$ all coincide (up to natural isomorphisms) for
$m,n$ sufficiantly large, and their common value will be called the
\underline{$C$-stabilization} of $\GA$, and denoted by $C^{\infty,\infty}\GA$.
\end{definition}
The proof of  the   following proposition is postponed until
sections 7--12:
\begin{prop}\label{prop:stable}
Let $\GA=(R,S)$ be a finite binary relation. Then
$C^{\infty,\infty}\GA$  is in a unique way the disjoint union
 of a finite number of sub-graphs of the following
types:
\begin{description}
\item[TYPE $C_N \;(N\geq 1)$] This  is an $N$-cycle consisting of
$N$ vertices $v_1,v_2,\cdots,v_N$,
 with relation consisting of the $N$ edges
 $$(v_1,v_2),(v_2,v_3),\cdots,(v_{N-1},v_N),(v_N,v_1)$$
$($If $N=1$ this consists of a single vertex $v$, together with a
loop $(v,v)$.$)$
\item[TYPE $ L_N \;(N\geq 1)$] This consists of $N$ vertices $v_1,v_2,\cdots,v_N$, with
relation consisting of the set of $N-1$ edges
$$(v_1,v_2),(v_2,v_3),\cdots,(v_{N-1},v_N)$$
$($If $N=1$ this consists of a single vertex $v$, with the set of
edges empty.$)$
\end{description}
\end{prop}
\smallskip

\par We are now ready to state  a second part of the main theorem Th.\ref{th:main}:
\setcounter{theorem}{3}
\begin{theorem}
$({\bf Main\; Theorem, Part\; Two})$ Let $\GA=(R\subseteq S\times
S)$ be a finite binary relation, with $K$-linearization
$K(\GA)$.With the notation of Prop.\ref{prop:stable}: Let
$C^{\infty,\infty}\GA$ be the disjoint union of $p$ graphs
$L(m_1),L(m_2),\cdots,L(m_p)$, together with $q$ cyclic graphs
$C(n_1),C(n_2),\cdots,C(n_q)$. Then:
\begin{description}
\item[a)]The indecomposable summands of Type $T_n$ which occur in the
Kronecker decomposition of $L(\GA)$ are given (multiplicity included) by the list
$$T_{m_1},T_{m_2},\cdots,T_{m_p}$$
---in other words, for every positive integer $n$, $t_n(K(\GA))$ equals
 the number of times $n$ is repeated in the sequence $m_1,m_2,\cdots,m_p$.
 \item[b)] The direct sum of those indecomposable summands of type $S(p(X)^n)$
 in the Kronecker decomposition of $K(\GA)$, is isomorphic to
 \begin{equation}\label{E:nonsing}
 \bigoplus_{i=1}^q S(X^{n_i}-1)
 \end{equation}
 (a fact which furnishes the summands of type $S(p(X)^n)$ via Eqn.\ref{E:CRT}.)
\end{description}
\end{theorem}

There is, however, one Kronecker multiplicity still undetermined by the first two
parts of Th.\ref{th:main}, namely the multiplicity $^0t^0_0(\GA)$ in $K(\GA)$
of the rather exceptional type
$$^0T^0_0=[0,0:K\to 0]\;.$$
(What makes this type unusual, is perhaps, that it is the only
indecomposable type in Kronecker's list $\mathcal{L}(K)$, which is
\emph{not} reduced, i.e. not a linear relation, in the sense
explained in Section \ref{S:7} below). Note that this type is
  not covered by equations (\ref{E:ZTZodd}) and  (\ref{E:ZTZeven})
 (which are explicitly stated to be restricted to $^0t^0_n$ with $n>0$).

 There will next be presented two different methods for the
 computation of this remaining Kronecker multiplicity
 $^0t^0_0(\GA)$, one of which uses all the other multiplicities of
 $\Gamma$. A nice check on the computation is then furnished by
 the requirement that these two methods yield the same result.

 The first of these methods is given by:
 \setcounter{theorem}{6}
 \begin{prop}\label{mult:ZTZ}
 Let
 $$\tau=[\mu,\nu:E\to V]$$
 be a transformation-pair over K; then
 $$ ^0t^0_0(\tau)=\dim(Ker\; \mu \cap Ker\; \nu)\;.$$
 \end{prop}

In order to explain the second method for computing  $^0t^0_0(\GA)$,
we first need the concept of \emph{edge-number} and \emph{vertex-number}:

\begin{definition}\label{def:edge}
Let $\tau$ be an indecomposable transformation-pair over $K$ in Kronecker's
list $\mathcal{L}(K)$; then the \underline{{\bf edge-number} $ \mathcal{E}(\tau)$}
and the \underline{{\bf vertex-number} $ \mathcal{V}(\tau)$}of $\tau$, are defined
as follows:
\begin{itemize}\label{edgeCounts}
\item For $n>0$,
$$\mathcal{E}(^0T_n)= \mathcal{V}(^0T_n)=n,\;\mathcal{E}(T^0_n)=
\mathcal{V}(T^0_n)=n$$
\item For $n\geq 0$,
$$\mathcal{E}(T_n)=n,\;\mathcal{V}(T_n)=n+1\;.$$
\item For $n\geq 0$;
$$ \mathcal{E}(^0T^0_n)=n+1,\mathcal{V}(^0T^0_n)=n$$
\item For $n>0$, and p a monic irreducible in $K[X]$ with $p\neq X$,
$$ \mathcal{E}(S(p^n))=n\cdot \mbox{ deg }p=\mathcal{V}(S(p^n))$$
\end{itemize}

Also, if $\GA=(R\subseteq S\times S)$ is a finite binary
relation, we define its  \underline{{\bf edge-number} $\mathcal{E}(\GA)$}
to be the cardinality of $R$, and its
\underline{{\bf vertex-number} $\mathcal{V}(\GA)$} to be the
cardinality of $S$:
$$\mathcal{E}(\GA)=\#(R)\mbox{ and }\mathcal{V}(\GA)=\#(S)\;.$$
\end{definition}

\setcounter{theorem}{3}
\begin{theorem}
$(${\bf Main Theorem, Part Three}$)$ Let $\GA=(R\subseteq S\times
S)$ be a finite binary relation, with $K$-linearization $K(\GA)$.
Let
$$K(\GA)=\bigoplus_{\tau\in\mathcal{L}(K)}[\GA:\tau]\tau$$
be the Kronecker decomposition of $K(\GA)$ with respect to Kronecker's list
$\mathcal{L}(K)$ of indecomposable $K$-transformation-pairs; then
\begin{equation}\label{E:edgesAdd}
\mathcal{E}(\GA)=\sum_{\tau\in\mathcal{L}(K)}[\GA:\tau]\mathcal{E}(\tau)\;,
\end{equation}
and
\begin{equation}\label{E:verticesAdd}
\mathcal{V}(\GA)=\sum_{\tau\in\mathcal{L}(K)}[\GA:\tau]\mathcal{V}(\tau)\;.
\end{equation}
\end{theorem}

\bigskip
\noindent {\bf NOTE:} By definition, $^0T^0_0$ has edge-number 1,
so we may solve eqn.(\ref{E:edgesAdd}) for $^0t^0_0(K\GA)$ in
terms of the other Kronecker multiplicities of $K\GA$ ( which we
may regard as already determined by Parts One and Two of
Th.\ref{th:main}), thus obtaining:
\begin{eqnarray}\label{solveForLast}
^0t^0_0&=&\mathcal{E}(\GA)-\sum_{n>0}n(^0t_n+t^0_n)-\sum_{n >
0}n\cdot t_n-\\
&&-\sum_{n>0}(n+1)\cdot^0\!t^0_n-\sum_{p,n>0}n\cdot \mbox{deg
}p\cdot[K(\GA):S(p^n)] \nonumber
\end{eqnarray}
with all multiplicities being evaluated for $K(\GA)$. Also, since $^0T^0_0$ has
vertex-number 0, we may also rewrite eqn.(\ref{E:verticesAdd}) in the form
\begin{eqnarray}\label{check}
\mathcal{V}(\GA)&=&\sum_{n>0}n(^0t_n+t^0_n)+\sum_{n\geq
0}(n+1)\cdot t_n+\\
&&+\sum_{n>0}n\cdot^0\!\!t^0_n+\sum_{p,n>0}n\cdot \mbox{deg
}p\cdot[K(\GA):S(p^n)] \nonumber
\end{eqnarray}
(again, with all multiplicities being evaluated for $K(\GA)$), and
in this form it provides a nice check for the multiplicities
computed using Parts One and Two of Th.\ref{th:main}.

\bigskip
Clearly, once proved, the three parts of Th.\ref{th:main}
completely solve the problem formulated in Sections 1 and 2. The
remainder of this paper (with the exception of the next section)
will be devoted to the proof of Th.\ref{th:main} (together with
the proofs of: Prop.\ref{ClCrCommute}, the assumptions underlying
Definition \ref{stable}, Prop. \ref{prop:stable}, and
Prop.\ref{mult:ZTZ}).

\smallskip
Let us note the following immediate consequence of Th.\ref{th:main}:

\smallskip
\setcounter{theorem}{7}
\begin{theorem}\label{completeSet} Let $\GA,\GA'$ be
two finite binary relations; then necessary and sufficient
that $\GA$ and $\GA'$ be linearly equivalent (i.e.
that the transformation-pairs $K(\GA)$ and $K(\GA')$ be isomorphic)
 is that, first,
$$\gamma_i,j(\GA)=\gamma_i,j(\GA')\;,$$
for all natural numbers $i,j$ with $|i-j|\leq 2$;\\
\noindent secondly, that $C^{\infty,\infty}(\GA)$ and $C^{\infty,\infty}(\GA')$ have the
same decompositions (in accordance with Prop.\ref{prop:stable} above) into
pieces of types $L_N$ and $C_N$  ;

\noindent and thirdly that $\GA$ and $\GA'$ have the same number
of edges. \noindent $($We may also express this by saying, that
$\{\G{i}{j}(\GA)\}$,the type of $C^{\infty,\infty}(\GA)$---i.e.,
the number of graphs $C(i),L(j)$ in its decomposition according to
Prop.\ref{prop:stable}--- and the natural number
$\mathcal{E}(\GA)$, form a complete set of linear equivalence
invariants for a binary relation $\GA$.$)$
\end{theorem}

\smallskip
  {\bf Note:} For given $\GA$, these linear invariants of $\GA$
\underline{{\bf are completely independent}} \underline{{\bf of
$K$}}.\footnote{Of course, if $L(N)$ is one of the graphs of type
$L$ occurring in the decomposition of $C^{\infty,\infty}(\GA)$, we
must factor $X^N-1$ over $K$ to obtain the precise Kronecker
decomposition of $K(\Gamma)$---but this factorization is not
needed to obtain a complete list of $K$-linear invariants of
$\GA$. }

\section{Some Illustrative Examples}\label{S:EX}

For the purpose of better understanding the algorithm explained in
the preceding section, we shall in this section apply it (without
proof) to a few specific examples. The remainder of the paper
after this section, then establishes the correctness of this
algorithm.

\bigskip
Note that the $\gamma_{m,n}$ which occur in
equations (\ref{E:ZT}), (\ref{E:TZ}), (\ref{E:ZTZodd}), (\ref{E:ZTZeven})
of Th.\ref{th:main}, all satisfy
\begin{equation}\label{E:suitable}
|m-n|\leq 2.
\end{equation}
We shall call a lattice-point $(m,n)$ \underline{{\bf suitable}} if it
satisfies (\ref{E:suitable}). It will be convenient to fit the various
suitable $(m,n)$ into a single picture as follows:
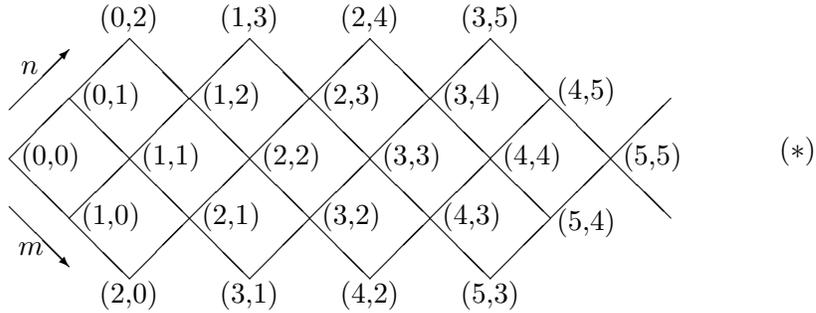
\setlength{\unitlength}{.8cm}
\begin{picture}(0,.8)(-1.7,2.3)
\put(0,0){\line(1,1){1 }} \put(0,0){\line(1,-1){1 }}
\put(1,1 ){\line(1,1){1 }} \put(1,1 ){\line(1,-1){1 }}
 \put(1,-1 ){\line(1,1){1 }} \put(1,-1 ){\line(1,-1){1 }}
\put(2,-2 ){\line(1,1){1 }} \put(2,2 ){\line(1,-1){1 }}
\put(2,0 ){\line(1,1){1 }} \put(2,0){\line(1,-1){1 }}
\put(3,-1 ){\line(1,1){1 }} \put(3,-1){\line(1,-1){1 }}
\put(3,1){\line(1,1){1 }} \put(3,1){\line(1,-1){1 }}
\put(4,-2){\line(1,1){1 }} \put(4,2){\line(1,-1){1 }}
 \put(4,0){\line(1,1){1 }} \put(4,0){\line(1,-1){1 }}
\put(5,1){\line(1,1){1 }} \put(5,1){\line(1,-1){1 }}
\put(5,-1){\line(1,1){1 }} \put(5,-1){\line(1,-1){1 }}
\put(6,-2){\line(1,1){1 }} \put(6,2){\line(1,-1){1 }}
\put(6,0){\line(1,1){1 }} \put(6,0){\line(1,-1){1 }}
\put(7,1){\line(1,1){1 }} \put(7,1){\line(1,-1){1 }}
\put(7,-1){\line(1,1){1 }} \put(7,-1){\line(1,-1){1 }}
\put(8,-2){\line(1,1){1 }} \put(8,2){\line(1,-1){1 }}
\put(8,0){\line(1,1){1 }} \put(8,0){\line(1,-1){1 }}
\put(9,-1){\line(1,1){1 }} \put(9,1){\line(1,-1){1 }}
\put(10,0){\line(1,1){1 }} \put(10,0){\line(1,-1){1 }}
\put(.2,-.1){(0,0)}\put(1.2,.9 ){(0,1)}\put(1.5,2.2){(0,2)}
\put(3.5,2,2){(1,3)}\put(5.5,2.2){(2,4)}\put(7.5,2.2){(3,5)}
\put(1.5,-2.4){(2,0)}\put(3.5,-2.4){(3,1)}\put(5.5,-2.4){(4,2)}\put(7.5,-2.4){(5,3)}
\put(1.2,-1.1 ){(1,0)}\put(2.2,-.1){(1,1)}\put(4.2,-.1){(2,2)}\put(6.2,-.1){(3,3)}
\put(8.2,-.1){(4,4)}\put(10.2,-.1){(5,5)}\put(3.2,-1.1){(2,1)}\put(5.2,-1.1){(3,2)}
\put(7.2,-1.1){(4,3)}\put(9.1,-1.2){(5,4)}
 \put(3.2,.9){(1,2)}\put(5.2,.9){(2,3)}\put(7.2,.9){(3,4)}\put(9.1,1){(4,5)}
\put(0,-.8){\vector(1,-1){1}}\put(.15,-1.6){$m$}
\put(0,.8){\vector(1,1){1}}\put(.2,1.4){$n$}
\put(12.8,0){($\ast$)}
\end{picture}

 \vspace{1.65in}
 For a given finite binary relation $\GA$, we obtain the
\underline{{\bf contraction-diagram }} $\mathcal{D}(\GA)$ for $\GA$ by inserting the integers
$\G{m}{n}$ for $(m,n)$ in the above diagram. For example, for the graph
$\GA_3$ of Example \ref{Ex4.4}, the pictures
in that Example show that
$$\G{0}{0}(\GA_3)=4,\G{1}{0}(\GA_3)=3=\G{0}{1}(\GA_3)$$
Also, it will readily be verified that
$$(C_l)^2\GA_3=C_l\GA_3,(C_r)^2\GA_3=C_r\GA_3$$
and that
$$C_lC_r\GA_3\cong \GA_2\cong C_rC_l\GA_3$$
where $\GA_2$ is the graph in Example 2.2; so we have
$$\G{2}{0}(\GA_3)=\G{0}{2}(\GA_3)=3,\G{1}{1}(\GA_3)=2\;.$$
Thus the contraction-diagram for $\GA_3$ is:
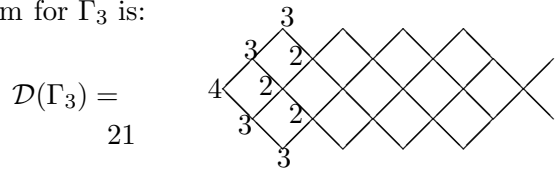
\setlength{\unitlength}{.4cm}
\begin{picture}(0,.8)(-1.7,2.3)
\put(-7,-.5){$\mathcal{D}(\GA_3)=$}
\put(0,0){\line(1,1){1 }} \put(0,0){\line(1,-1){1 }}
\put(1,1 ){\line(1,1){1 }} \put(1,1 ){\line(1,-1){1 }}
 \put(1,-1 ){\line(1,1){1 }} \put(1,-1 ){\line(1,-1){1 }}
\put(2,-2 ){\line(1,1){1 }} \put(2,2 ){\line(1,-1){1 }}
\put(2,0 ){\line(1,1){1 }} \put(2,0){\line(1,-1){1 }}
\put(3,-1 ){\line(1,1){1 }} \put(3,-1){\line(1,-1){1 }}
\put(3,1){\line(1,1){1 }} \put(3,1){\line(1,-1){1 }}
\put(4,-2){\line(1,1){1 }} \put(4,2){\line(1,-1){1 }}
 \put(4,0){\line(1,1){1 }} \put(4,0){\line(1,-1){1 }}
\put(5,1){\line(1,1){1 }} \put(5,1){\line(1,-1){1 }}
\put(5,-1){\line(1,1){1 }} \put(5,-1){\line(1,-1){1 }}
\put(6,-2){\line(1,1){1 }} \put(6,2){\line(1,-1){1 }}
\put(6,0){\line(1,1){1 }} \put(6,0){\line(1,-1){1 }}
\put(7,1){\line(1,1){1 }} \put(7,1){\line(1,-1){1 }}
\put(7,-1){\line(1,1){1 }} \put(7,-1){\line(1,-1){1 }}
\put(8,-2){\line(1,1){1 }} \put(8,2){\line(1,-1){1 }}
\put(8,0){\line(1,1){1 }} \put(8,0){\line(1,-1){1 }}
\put(9,-1){\line(1,1){1 }} \put(9,1){\line(1,-1){1 }}
\put(10,0){\line(1,1){1 }} \put(10,0){\line(1,-1){1 }}
\put(-.5,-.3){4}\put(.7,1){3}\put(1.9,2.1){3}
\put(.5,-1.5){3}\put(1.8,-2.6){3}
\put(1.2,-.2){2}\put(2.2,.8){2}\put(2.2,-1.15){2}
\end{picture}

\vspace{.9in} \noindent where all the intersections with no
integer marked, are assigned to 2. The reader may wish to verify,
as an exercise, the following somewhat abbreviated
contraction-diagrams for the graphs $\GA_1$ and $\GA_2$ in
Examples 2.1 and 2.2:

\begin{center}
\setlength{\unitlength}{.4cm}
\begin{picture}(0,0)(6,1)
\put(0,0){\line(1,1){1 }} \put(0,0){\line(1,-1){1}}
\put(1,1){\line(1,-1){1}} \put(1,-1){\line(1,1){1}}
\put(1,-1){\line(1,-1){1}} \put(2,0){\line(1,-1){1}}
\put(2,-2){\line(1,1){1 }} \put(-.6,-.6){3} \put(.3,-1.5){2}
\put(1.3,-2.5){2} \put(1.1,1.1){2} \put(2.1,.1){1}
 \put(3.1,-.9){1} \put(-4.5,0){$\mathcal{D}(\GA_1)=$}
 \put(4,0){,} \put(9.5,0){\line(1,1){1 }}
\put(9.5,0){\line(1,-1){1 }}
  \put(10.5,1 ){\line(1,-1){1 }}
\put(10.5,-1 ){\line(1,1){1 }}
 \put(10.1,1.3){1}
\put(10.1,-1.8){1}
\put(11.5,-.2){1}
\put(5,0){$\mathcal{D}(\GA_2)=$}
\put(9,-.3){2}
\end{picture}
\end{center}

\vspace{.4in}
Unfortunately, the three graphs $\GA_1$,$\GA_2$,$\GA_3$ studied up to this
point, seem too small to show all the features to be illustrated. It's time
for a slightly heftier example:
Consider the binary relation
$$\GA_4=(R_4,S_4)$$
on the set $S_4$ of integers from 0 to 17, represented by the following diagram:

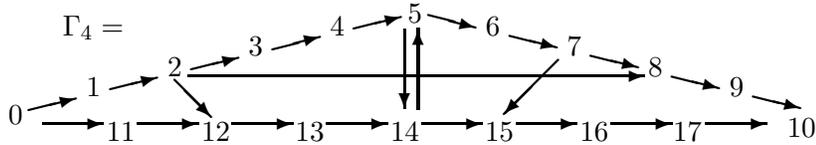
\setlength{\unitlength}{.9cm}
 \begin{picture}(4,5)(0,-3)
 \thicklines
\multiput(0,0)(1.2,.3){5}{\vector(4,1){.7}}
\put(5.6,1.3){5}
\multiput(5.9,1.4)(1.2,-.3){5}{\vector(4,-1){.7}}
\put(-.3,-.2){0}\put(11.2,-.4){10}
\multiput(.2,-.2)(1.4,0){8}{\vector(1,0){.9 }}
\put(.85,.2){1}\put(2.05,.5){2}\put(3.25,.8){3}
\put(4.45,1.1){4}\put(6.75,1.1){6}\put(7.95,.8){7}
\put(9.15,.5){8}\put(10.35,.2){9}
\put(1.15,-.45){11}\put(2.55,-.45){12}\put(3.94 ,-.45){13}\put(5.34,-.45){14}
\put(6.74,-.45){15}  \put(8.14,-.45){16}\put(9.52,-.45){17}
\put(5.55,1.2){\vector(0,-1){1.18}}\put(5.75,0){\vector(0,1){1.2}}
\put(2.15,.45){\vector(1,-1){.55}}\put(2.35,.5){\vector(1,0){6.75}}
\put(7.82,.74){\vector(-1,-1){.8}}\put(.5,1.1){$\GA_4=$}
\end{picture}

 \vspace{-.6in}
Let us next compute the left contraction
$$C_l\GA_4=^1\!\!(\GA_4)^0$$
 of this binary relation:

\noindent We first list all cases where two arrows emerge \emph{from} the
same vertex: they are furnished by
$$1\leftarrow 0\rightarrow11,\;6\leftarrow5\rightarrow14,\;8\leftarrow7\rightarrow15,
\;5\leftarrow14 \rightarrow15$$
together with the 3 arrows emerging from 2 to 3, 8, 12.

\noindent Thus the equivalence relation $\simL=^1\sim^0$ is here generated by
$$1\simL11,\;6\simL 14,\;8\simL15,\;5\simL15\mbox{ and }3\simL8\simL12\simL3$$
yielding the 12 $\simL$-equivalence-classes
$$\{0\},\;\{1,11\},\;\{2\},\;\{3,5,8,12,15\},\;\{4\},\;\{6,14\},\;\{7\},\;
\{9\},\;\{10\},\;\{13\},\;\{16\},\;\{17\},\;.$$
which make up the vertices of the left-contracted binary relation $C_l\GA_4$. Hence
 $\G{1}{0}(\GA_4)=12\;.$
In order to iterate this contraction process, and so compute the higher
$\gamma$'s, we
must also compute the contracted relation
$C_lR_4$ on these 12 equivalence classes; using Def.\ref{def:quotient} we obtain
the following diagram:

\setlength{\unitlength}{.7cm}
 \begin{picture}(4,5)(0,-1.4)
 \thicklines
\put(.55,0){\vector(1,0){.9}}\put(-.16,-.1){$\{0\}$}
\put(1.5,-.1){$\{1,11\}$}\put(1.8,.4){\vector(0,1){.75}}
\put(1.3,1.35){$\{2\}$}\put(2.3,.35){\vector(3,1){2.3}}
  \put(3.4,1.3){$\{3,5,8,12,15\}$} \put(2.15,1.45){\vector(1,0){1.2}}
  \put(5.3,1.2){\vector(-1,-1){.75}}\put(3.6 ,-.1){$\{13\}$}
  \put(4.7,0){\vector(1,0){.6}}\put(5.3,-.1){$\{6,14\}$}
  \put(5.9,1.2){\vector(0,-1){.8 }}\put(6.2,.4){\vector(0,1){.8 }}
 \put(6.93,-.05){\vector(1,0){.5}}
\put(7.38,-.1){$\{7\}$}\put(7.4,.3){\vector(-1,1){.77}}
\put(7,1.3){\vector(3,-1){3}}\put(10,-.1){$\{9\}$}
\put(10.88,.05){\vector(1,0){1}}\put(11.9,-.1){$\{10\}$}
\put(11.9,1.3){$\{17\}$}\put(12.4,1.2){\vector(0,-1){.8}}
\put(6.9,1.5){\vector(1,0){2.6}}\put(9.5,1.3){$\{16\}$}
\put(10.65,1.5){\vector(1,0){1.2}}\put(5.4 ,1.75){\vector(0,1){.8}}
\put(5.1,2.75){$\{4\}$}\put(5.7,2.55){\vector(0,-1){.8}}
\put(.6,2.75){$C_l \GA_4=:$}

\end{picture}
A similar straightforward process can be used to compute the right contraction
$C_r\GA_4$. Here, we must begin
by listing all cases where two arrows converge \emph{to} the same vertex
---for $\GA_4$, one such example is
$$11\rightarrow12\leftarrow 2,\mbox{ whence }11\simR2\;.$$
The integers from 0 to 17 then divide into 12 equivalence-classes with respect to
$\simR(\GA_4)$, namely
$$\{0\},\{1\},\{2,4,7,11,14\},\{3\},\{5,13\},\{6\},\{8\},\{9,17\},\{10\},\{12\},\{15\},\{16\}$$
so that $\G{0}{1}(\GA)=12$.

Let us continue this process, computing the binary relation
$$C_l^LC_r^M(\GA_4)=(R_4^{L,M},S^{L,M})$$
for the first few $(L,M)$ which are `suitable', i.e., for which $|L-M|\leq 2$.
We obtain a steadily coarsening collection of partitionings $S^{L,M}$of
$$S_4=\{0,1,\cdots,17\}$$
of which the first three are given above, while also \\
$S^{2,0}=\{\{0\},\{1,11\},\{2,3,5,7,8,12,15\},\{4,6,9,13,14,16\}
,\{10\},\{17\}\}  $ \\
$S^{1,1}=\{\{0\},\{1,2,4,6,7,11,14\},\{3,5,8,12,13,15\},
\{9,17\},\{10\},\{16\} \} $ \\
 $S^{0,2,}=\{\{0,1,3,5,6,13\},\{2,4,7,11,12,14\},
 \{8,16\},\{9,17\},\{10\},\{15\}   \} $ \\
$S^{2,1}=\{S_4\backslash \{0,10\},\{0\},\{10\}   \} $ \\
$S^{1,2}=S^{1,3}=\{S_4\backslash \{9,10,17\},\{9,17\},\{10\}   \} $ \\
$S^{2,2}=S^{2,3}=S^{2,4}=\{S_4\backslash \{10\},\{10\}   \} $ \\
$S^{3,1,}=\{S_4\backslash \{0\},\{0\}   \} $ \\
---and where all further $S^{L,M}$ with $(L,M)$ suitable
consist of the single equivalence class
$$S_4=\{0,1,\cdots,17\}$$.

\medskip
Thus, we obtain the  values of $\gamma_{m,n}(\GA_4)$ in
 the following
contraction-diagram for $\GA_4$:

\setlength{\unitlength}{.8cm}
\begin{picture}(0,.8)(-2.5,2.3)
\put(0,0){\line(1,1){1 }} \put(0,0){\line(1,-1){1 }}
\put(1,1 ){\line(1,1){1 }} \put(1,1 ){\line(1,-1){1 }}
 \put(1,-1 ){\line(1,1){1 }} \put(1,-1 ){\line(1,-1){1 }}
\put(2,-2 ){\line(1,1){1 }} \put(2,2 ){\line(1,-1){1 }}
\put(2,0 ){\line(1,1){1 }} \put(2,0){\line(1,-1){1 }}
\put(3,-1 ){\line(1,1){1 }} \put(3,-1){\line(1,-1){1 }}
\put(3,1){\line(1,1){1 }} \put(3,1){\line(1,-1){1 }}
\put(4,-2){\line(1,1){1 }} \put(4,2){\line(1,-1){1 }}
 \put(4,0){\line(1,1){1 }} \put(4,0){\line(1,-1){1 }}
\put(.24,-.1){18}\put(1.15,.9){12}\put(1.15,-1.1){12}
\put(2.25,-.1){6}\put(3.2,.9){3}\put(3.2,-1.1){3}
\put(4.24,-.1 ){2}\put(1.85,2.1){6}\put(1.85,-2.4){6}
\put(-2.5,0){$\mathcal{D}(\GA_4)=$}
\put(5,1){\line(1,1){1 }} \put(5,1){\line(1,-1){1 }}
\put(5,-1){\line(1,1){1 }} \put(5,-1){\line(1,-1){1 }}
\put(6,-2){\line(1,1){1 }} \put(6,2){\line(1,-1){1 }}
\put(6,0){\line(1,1){1 }} \put(6,0){\line(1,-1){1 }}
\put(7,1){\line(1,1){1 }} \put(7,1){\line(1,-1){1 }}
\put(7,-1){\line(1,1){1 }} \put(7,-1){\line(1,-1){1 }}
\put(8,-2){\line(1,1){1 }} \put(8,2){\line(1,-1){1 }}
\put(8,0){\line(1,1){1 }} \put(8,0){\line(1,-1){1 }}
\put(9,-1){\line(1,1){1 }} \put(9,1){\line(1,-1){1 }}
\put(10,0){\line(1,1){1 }} \put(10,0){\line(1,-1){1 }}
\put(3.85,2.1){3}\put(3.9,-2.4){2}\put(5.18,-1.15){1}
\put(5.25,.9){2}\put(6.18,-.15){1}\put(5.85,2.1){2}
\put(7.2,.86){1}
 \put(0,-.8){\vector(1,-1){1}}\put(.15,-1.6){$m$}
\put(0,.8){\vector(1,1){1}}\put(.2,1.4){$n$}
\end{picture}

\vspace{1.6in} \noindent (abbreviated by the same convention as
before, whereby all `suitable' vertices not shown or not labelled,
are understood to have assigned as labels the minimal
$\gamma$-value, which here is 1).

Another useful further result of this computation, is that ---as
a special case of Prop.\ref{prop:stable}--- the graphs
$C_l^LC_m^M(\GA_4)$ repeatedly contracted from $\GA_4$ eventually
stablize
 at the graph
$C^{\infty,\infty}(\GA_4)=C_1$, consisting of one vertex, and one
loop at that vertex.

Thus,we have computed $\mathcal{D}(\GA_4)$ and $C^{\infty,\infty}(\GA_4)$; it only
remains to note the number
$$\mathcal{E}(\GA_4)=23$$
of edges in $\GA_4$, to have obtained a complete set of linear
equivalence invariants, according to Theorem \ref{completeSet}.

\bigskip
\noindent {\bf NOTE:} A slight strengthening of Prop
\ref{ClCrCommute} furnishes a useful repeated check during such
computations, which works as follows: Consider, for example, the
computation of
$$\G{2}{2}(\GA_4)=\mbox{ number of   vertices in }C_l^2C_r^2(\GA_4)\;.$$
The check in question consists in noting that we have two
ways to construct the `multi-contracted' graph in question,
namely as $C_l(C_lC_r^2)$ and as $C_r(C_l^2C_r)$, and also in
noting that, if as above we compute these graphs in terms of partitionings
of $S_4$, then we have not only graph-isomorphism between them
 (as asserted by Prop.\ref{ClCrCommute}), but even more:
actual \emph{identity} of the partitionings involved (whence the binary
relations involved, as determined by Def.\ref{def:quotient}, also coincide).
This stronger version of Prop \ref{ClCrCommute} is actually what is proved below, in
\S\ref{S:6}.

\bigskip
Having thus obtained a complete set of linear equivalence
invariants for $\GA_4$, the assertions of the preceding \S \ tell
us how to compute from them, the Kronecker invariants for the
$K$-transformation pair $K(\GA_4)$. Let us now follow the
instructions for doing so:

\bigskip
We first direct our attention to Part One of Th.\ref{th:main}, which expresses
\emph{some} of  the Kronecker invariants ---namely
$$^0t_n(K(\GA)),t^0_n(K(\GA)),^0\!t^0_n(K(\GA))\mbox{  (all with  }n>0)$$
---in terms of the $\G{m}{n}(\GA)$, ie of the contraction-diagram for $\GA$. Equations
 (\ref{E:ZT})-(\ref{E:ZTZeven}) of this theorem, may perhaps be more easily visualized in term
of the contraction-diagram, as follows.

At this point, the reader is asked to examine carefully the picture $(\ast)$
at the beginning of this section. We shall see that this picture
may usefully be sub-divided into
smaller parts according to three different methods. We begin with the simplest
method of sub-division:

\noindent {\bf SUB-DIVISION ONE:}\\
This sub-divides the picture $(\ast)$ into
the coordinate squares of the $(m,n)$ lattice---as follows:
\vspace{1in}

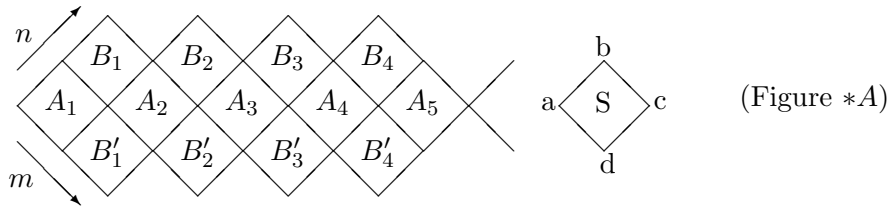
\setlength{\unitlength}{.6cm}
\begin{picture}(0,.9)(-1.7,2.3)
\put(0,0){\line(1,1){1 }} \put(0,0){\line(1,-1){1 }}
\put(1,1 ){\line(1,1){1 }} \put(1,1 ){\line(1,-1){1 }}
 \put(1,-1 ){\line(1,1){1 }} \put(1,-1 ){\line(1,-1){1 }}
\put(2,-2 ){\line(1,1){1 }} \put(2,2 ){\line(1,-1){1 }}
\put(2,0 ){\line(1,1){1 }} \put(2,0){\line(1,-1){1 }}
\put(3,-1 ){\line(1,1){1 }} \put(3,-1){\line(1,-1){1 }}
\put(3,1){\line(1,1){1 }} \put(3,1){\line(1,-1){1 }}
\put(4,-2){\line(1,1){1 }} \put(4,2){\line(1,-1){1 }}
 \put(4,0){\line(1,1){1 }} \put(4,0){\line(1,-1){1 }}
\put(5,1){\line(1,1){1 }} \put(5,1){\line(1,-1){1 }}
\put(5,-1){\line(1,1){1 }} \put(5,-1){\line(1,-1){1 }}
\put(6,-2){\line(1,1){1 }} \put(6,2){\line(1,-1){1 }}
\put(6,0){\line(1,1){1 }} \put(6,0){\line(1,-1){1 }}
\put(7,1){\line(1,1){1 }} \put(7,1){\line(1,-1){1 }}
\put(7,-1){\line(1,1){1 }} \put(7,-1){\line(1,-1){1 }}
\put(8,-2){\line(1,1){1 }} \put(8,2){\line(1,-1){1 }}
\put(8,0){\line(1,1){1 }} \put(8,0){\line(1,-1){1 }}
\put(9,-1){\line(1,1){1 }} \put(9,1){\line(1,-1){1 }}
\put(10,0){\line(1,1){1 }} \put(10,0){\line(1,-1){1 }}
\put(0,-.8){\vector(1,-1){1.4}}\put(-.2,-1.8){$m$}
\put(0,.8){\vector(1,1){1.4}}\put(-.05,1.4){$n$}
\put(.6,-.1){$A_1$}\put(2.6,-.1){$A_2$}\put(4.6,-.1){$A_3$}
 \put(6.6,-.1){$A_4$}\put(8.6,-.1){$A_5$}
\put(1.6,.9){$B_1$}\put(3.6 ,.9){$B_2$}\put(5.6 ,.9){$B_3$}\put(7.6,.9){$B_4$}
\put(1.6,-1.2){$B_1'$}\put(3.6,-1.2){$B_2'$}\put(5.6,-1.2){$B_3'$}
\put(7.6,-1.2){$B_4'$} \put(12,0){\line(1,1){1}}\put(12,0){\line(1,-1){1}}
\put(13,1){\line(1,-1){1}}\put(13,-1){\line(1,1){1}}
\put(12.8,-.15){S}
\put(11.6,-.1){a}\put(12.8,1.1){b}\put(14.1,-.08){c}
\put(12.9,-1.5){d}\put(16,0){(Figure $\ast A$)}
\end{picture}

\vspace{1.2in}
If $S$ is one of these squares, with vertices
$$a=(a_1,a_2),b=(b_1,b_2),c=(c_1,c_2),d=(d_1,d_2)$$ as indicated
in the above figure, let us define the \underline{{\bf $\GA$-content}} of $S$ to be
$$|S|_{\GA}:=\gamma_a(\GA)-\gamma_b(\GA)-\gamma_d(\GA)+\gamma_c(\GA);.$$
For instance, the $\GA$-content of the square $A_1$ is
$$|A_1|_{\GA}=\G{0}{0}(\GA)-\G{1}{0}(\GA)-\G{0}{1}(\GA)+\G{1}{1}(\GA)$$
which by equation (\ref{E:ZTZodd}) is equal to the Kronecker
invariant $^0t^0_1(K\GA)$. We may similarly visualize the
remaining information in equation (\ref{E:ZTZodd}): which tells
us that the $\GA$-contents of the central squares
$A_1$,$A_2$,$A_3$,etc. in the preceding Figure$\ast A$, are equal
respectively to the Kronecker invariants
$$^0t^0_1(K\GA),\;^0t^0_3(K\GA),\;^0t^0_5(K\GA),\cdots$$
with \emph{odd} subscripts.

\medskip Applying this method to the particular case $\GA_4$, we see that
values for $^0t^0_{2p+1}$ (\emph{odd} subscripts) are furnished by the
central squares $A_1, A_2, A_3,\cdots$ in $\mathcal{D}(\GA_4)$, as follows:
$$^0t^0_1(K\GA_4)=|A_1|_{K\GA_4}=18+6-12-12=0,\;^0t^0_3(K\GA_4)=6-3-3+2=2$$
while
$$^0t^0_5(K\GA_4)=2-2-1+1=0,\;^0t^0_7(K\GA_4)=^0\!\!t^0_9(K\GA_4)=\cdots=1-1-1+1=0$$

Similarly, applying these facts to the abbreviated contraction-diagrams computed
above for the binary relations $\GA_1$,$\GA_2$,$\GA_3$, we obtain \\
$ ^0t^0_1(K\GA_1)=|A_1|_{K\GA_1}=3-2-1+1=1,^0\!t^0_1(K\GA_2)=2-1-1+1=1,$
\\
$^0\!t^0_1(K\GA_3)=4-3-3+2=0$.\\
This is in agreement with the fact (obtained by direct \emph{ad
hoc} methods in \S2) that each of $K(\GA_1),K(\GA_2)$ has in its
Kronecker decomposition, exactly one direct summand of type
$^0T^0_1$.

\medskip
We may also use Figure $\ast A$ (together with the concept of
$\GA$-content) to visualize equation (\ref{E:ZTZeven}). In these
terms, this equation asserts two things:

\noindent A) The non-central off-diagonal squares in Figure $\ast A$ occur in pairs
$$B_1 \mbox{ and } B_1';\;\;B_2 \mbox{ and }B_2';\;\;B_3 \mbox{ and }B_3';\;\cdots$$
for which both elements in a pair have the same $\GA$-content:
$$|B_n|_{\GA}=|B_n'|_{\GA}\mbox{ for all }n\geq 1\;.$$

\noindent B) These common $\GA$-contents, are equal respectively to the
Kronecker invariants
$$^0t^0_2(K\GA),\;^0t^0_4(K\GA),\;^0t^0_6(K\GA),\cdots$$
with \emph{even} subscripts.

\smallskip For example, here is what we get if we apply these
assertions to the contraction-diagram for $\GA_4$:
$$|B_1|_{\GA_4}=12-6-6+3=|B_1'|_{\GA_4}=3=^0\!\!t^0_2(K\GA_4)$$
and all $^0t^0_{2p}(K\GA_4)$ with $p>1$ are 0--- e.g.,
$$|B_2|_{\GA_4}=3-3-2+2=|B_2'|_{K\GA_4}=3-2-2+1=0=^0\!\!t^0_4(K\GA_4)$$

So much for equations (\ref{E:ZTZodd}) and (\ref{E:ZTZeven}). In
order to find a similar interpretation of equation (\ref{E:ZT}),
which furnishes the values of the multiplicities $^0t_n$, it is
necessary to sub-divide a portion of the picture ($\ast$) at the
beginning of this section, into quadrilaterals by a second method.
Here, unlike Sub-Division One, where the entire diagram ($\ast$)
was covered by the quadrilaterals $A_i$,$B_i$,$B_i'$, we are led
by equation (\ref{E:ZT}) to introduce quadrilaterals
$$C_1,\;C_1',\;C_2,\;C_2',\;C_3,\;C_3',\cdots$$
which cover roughly $2/3$ of ($\ast$), as follows. The
quadrilaterals $C_1,C_1'$;\ $C_2,C_2';$\ $C_3,C_3';\;\cdots$ in
this figure are all parallelograms, with sides either vertical or
parallel to the indicated $n$-axis. (In addition to these
quadrilaterals, the picture also contains a number of isosceles
right triangles, which seem irrelevant to our present purposes.)

{\bf SUB-DIVISION TWO:}

 \setlength{\unitlength}{.9cm}
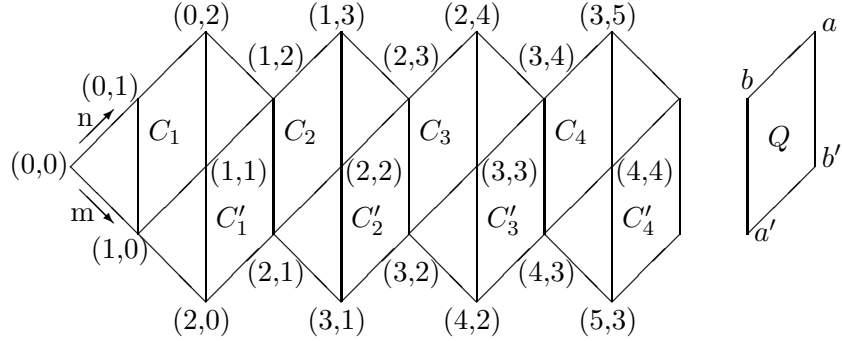
\begin{picture}(0,.9)(-1,1.7)
 \put(0,0){\line(1,1){1 }} \put(0,0){\line(1,-1){1 }}
\put(1,1){\line(0,-1){2}}\put(1,1){\line(0,-1){2}}
\put(1,1){\line(1,1){1}}\put(1,-1){\line(1,1){1}}\put(2,2){\line(0,-1){2}}
 \put(2,0){\line(1,1){1}}\put(2,0){\line(0,-1){2}}
 \put(3,1){\line(0,-1){2}}\put(2,-2){\line(1,1){1}}
\put(3,-1){\line(0,1){2}}\put(3,-1){\line(1,1){1}}
\put(1.15 ,.4){$C_1$}\put(2.15,-.9){$C_1'$}
\put(3,1){\line(1,1){1}}\put(4,2){\line(0,-1){2}}
\put(3.15,.4){$C_2$}
\put(4,-2){\line(0,1){2}}\put(4,-2){\line(1,1){1}}\put(5,-1){\line(0,1){2}}
\put(4,0){\line(1,1){1}}\put(4.15,-.9){$C_2'$}
\put(5,1){\line(1,1){1}}\put(5,-1){\line(1,1){1}}\put(6,0){\line(0,1){2}}
\put(5.15,.4){$C_3$}\put(6,0){\line(0,-1){2}}\put(6,0){\line(1,1){1}}
\put(6,-2){\line(1,1){1}}\put(7,-1){\line(0,1){2}}\put(6.15,-.9){$C_3'$}
\put(7,-1){\line(1,1){1}}\put(7,1){\line(1,1){1}}\put(8,0){\line(0,1){2}}
\put(7.15,.4){$C_4$}\put(8,-2){\line(0,1){2}}\put(8,-2){\line(1,1){1}}
\put(8,0){\line(1,1){1}}\put(9,-1){\line(0,1){2}}\put(8.15,-.9){$C_4'$}
\put(-.9,-.1){(0,0)}\put(.15,1.05 ){(0,1)}\put(1.5,2.1){(0,2)}
\put(3.5,2.1){(1,3)}\put(5.5,2.1){(2,4)}\put(7.5,2.1){(3,5)}
\put(1.5,-2.4){(2,0)}\put(3.5,-2.4){(3,1)}\put(5.5,-2.4){(4,2)}\put(7.5,-2.4){(5,3)}
\put(.3,-1.35 ){(1,0)}\put(2.05,-.2){(1,1)}\put(4.05,-.2){(2,2)}\put(6.05,-.2){(3,3)}
\put(8.05,-.2){(4,4)}\put(2.6,-1.7){(2,1)}\put(4.6,-1.7){(3,2)}
\put(6.6,-1.7){(4,3)} \put(2.6,1.5){(1,2)}
\put(4.6,1.5){(2,3)}\put(6.6,1.5){(3,4)}

\put(2,2){\line(1,-1){1}}\put(1,-1){\line(1,-1){1}}
\put(4,2){\line(1,-1){1}}\put(6,2){\line(1,-1){1}}\put(8,2){\line(1,-1){1}}
\put(3,-1){\line(1,-1){1}}\put(5,-1){\line(1,-1){1}}\put(7,-1){\line(1,-1){1}}
\put(10,-1){\line(0,1){2}}\put(10,-1){\line(1,1){1}}
\put(10,1){\line(1,1){1}}\put(11,0){\line(0,1){2}}
\put(10.3,.3){$Q$}
\put(10.1,-1.1){$a'$}\put(9.9,1.15){$b$}
\put(11.1,2){$a$}\put(11.1,0){$b'$}
\put(.14,-.34){\vector(1,-1){.5}}\put(0,-.8){m}
\put(.14,.34){\vector(1,1){.5}}\put(.1,.6){n}
\end{picture}

\vspace{1.8in}

 If $Q$ is one of these quadrilaterals, with vertices
$$a=(a_1,a_2),\;a'=(a_1',a_2'),\;b=(b_1,b_2),\;b'=(b_1',b'_2)$$
as indicated in the above figure, let us define the \underline{{\bf $\Gamma$-content of}}
$Q$ to be
$$|Q|_{\Gamma}:=\gamma_a(\GA)+\gamma_{a'}(\GA)-\gamma_b(\GA)-\gamma_{b'}(\GA)\;,$$
For instance, the $\GA$-content of $C_1$ is
$$|C_1|_{\GA}=\gamma_{0,2}(\GA)+\gamma_{1,0}(\GA)-
\gamma_{0,1}(\GA)-\gamma_{11}(\GA)$$
while that of $C_1'$ is
$$|C_1'|_{\GA}=\gamma_{1,2}(\GA)+\gamma_{2,0}(\GA)-\gamma_{1,1}(\GA)-\gamma_{2,1}(\GA)$$
and equation (\ref{E:ZT}) (for $n=1$) asserts that these two
expressions are equal, and that their common value is the
Kronecker multiplicity $^0t_1(K\GA)$. We may similarly visualize
the remaining information in equation  (\ref{E:ZT}), which asserts that,
for all positive $n$, $C_n$ and $C_n'$ have the same
$\GA$-content, and that this equals the Kronecker multiplicity
$^0t_n$ of $^0T_n$ in $K\GA$:
$$|C_n|_{\GA}=|C_n'|_{\GA}=^0\!\!t_n(K\GA)$$
For the particular case $\GA_4$, these facts become:
$$|C_1|_{\GA}=5+12-12-6=0=3+6-6-3=|C_1'|$$
giving the value 0 for $^0t_1(K\GA_4)$; and similarly
$$^0t_2(K\GA_4)=1,\;^0t_n(K\GA_4)=0\mbox{ for }n>2$$

The final information remaining to be discussed in
Th.\ref{th:main} is that furnished by equation (\ref{E:TZ}), which
we shall visualize in terms of the following

\noindent {\bf SUB-DIVISION THREE:}

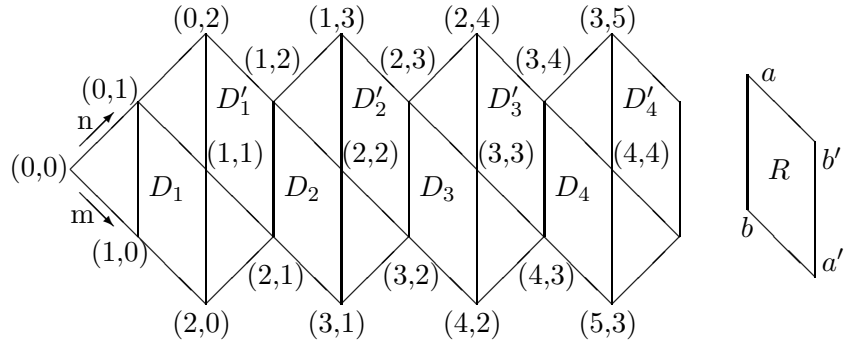
\setlength{\unitlength}{.9cm}
\begin{picture}(0,.9)(-1,2.3)
 \put(0,0){\line(1,1){1 }} \put(0,0){\line(1,-1){1 }}
\put(1,1){\line(0,-1){2}}
\put(1,-1){\line(1,-1){1}}\put(1,1){\line(1,-1){1}}\put(2,-2){\line(0,1){2}}
 \put(2,0){\line(1,-1){1}}\put(2,0){\line(0,1){2}}
 \put(3,-1){\line(0,1){2}}\put(2,2){\line(1,-1){1}}
\put(3,1){\line(0,-1){2}}\put(3,1){\line(1,-1){1}}
\put(1.15 ,-.4){$D_1$}\put(2.15,.9){$D_1'$}
\put(3,-1){\line(1,-1){1}}\put(4,-2){\line(0,1){2}}
\put(3.15,-.4){$D_2$}
\put(4,2){\line(0,-1){2}}\put(4,2){\line(1,-1){1}}\put(5,1){\line(0,-1){2}}
\put(4,0){\line(1,-1){1}}\put(4.15,.9){$D_2'$}
\put(5,-1){\line(1,-1){1}}\put(5,1){\line(1,-1){1}}\put(6,0){\line(0,-1){2}}
\put(5.15,-.4){$D_3$}\put(6,0){\line(0,1){2}}\put(6,0){\line(1,-1){1}}
\put(6,2){\line(1,-1){1}}\put(7,1){\line(0,-1){2}}\put(6.15,.9){$D_3'$}
\put(7,1){\line(1,-1){1}}\put(7,-1){\line(1,-1){1}}\put(8,0){\line(0,-1){2}}
\put(7.15,-.4){$D_4$}\put(8,2){\line(0,-1){2}}\put(8,2){\line(1,-1){1}}
\put(8,0){\line(1,-1){1}}\put(9,1){\line(0,-1){2}}\put(8.15,.9){$D_4'$}
\put(-.9,-.1){(0,0)}\put(.15,1.05 ){(0,1)}\put(1.5,2.1){(0,2)}
\put(3.5,2.1){(1,3)}\put(5.5,2.1){(2,4)}\put(7.5,2.1){(3,5)}
\put(1.5,-2.4){(2,0)}\put(3.5,-2.4){(3,1)}\put(5.5,-2.4){(4,2)}\put(7.5,-2.4){(5,3)}
\put(.3,-1.35 ){(1,0)}\put(2,.1){(1,1)}\put(4 ,.1){(2,2)}\put(6 ,.1){(3,3)}
\put(8 ,.1){(4,4)}\put(2.6,-1.7){(2,1)}\put(4.6,-1.7){(3,2)}
\put(6.6,-1.7){(4,3)} \put(2.55,1.5){(1,2)}
\put(4.55,1.5){(2,3)}\put(6.55,1.5){(3,4)}

\put(2,-2){\line(1,1){1}}\put(1,1){\line(1,1){1}}
\put(4,-2){\line(1,1){1}}\put(6,-2){\line(1,1){1}}\put(8,-2){\line(1,1){1}}
\put(3,1){\line(1,1){1}}\put(5,1){\line(1,1){1}}\put(7,1){\line(1,1){1}}
\put(10,1.4){\line(0,-1){2}}\put(10,1.4){\line(1,-1){1}}
\put(10,-.6){\line(1,-1){1}}\put(11,.4){\line(0,-1){2}}
\put(10.3,-.15){$R$}
\put(10.2,1.3){$a$}\put(9.9,-1 ){$b$}
\put(11.1,-1.5){$a'$}\put(11.1,0){$b'$}
\put(.14,-.34){\vector(1,-1){.5}}\put(0,-.8){m}
\put(.14,.34){\vector(1,1){.5}}\put(.1,.6){n}
\end{picture}

\vspace{1.8in}
Here, we have the diagram $(\ast)$ sub-divided into isosceles right triangles,
together with a collection of paired parallelograms
$$D_1\mbox{ and }D_1';\;D_2\mbox{ and }D_2';\;D_3\mbox{ and }D_3'\mbox{; etc.}$$
each of whose sides is either vertical, or parallel to the indicated m-axis. If R is
one of these parallelograms, with vertices as indicated, we define the
\underline{{\bf $\Gamma$-content of}}
$R$ to be
$$|R|_{\Gamma}:=\gamma_a(\GA)+\gamma_{a'}(\GA)-\gamma_b(\GA)-\gamma_{b'}(\GA)\;,$$
It is readily verified that equation (\ref{E:TZ}) is equivalent to
asserting that the two following statements hold for all positive
integers $n$ and all binary relations $\GA$:
\begin{description}
\item{Firstly,} the two paired parallelograms $D_n$ and $D_n'$ of this third
subdivision, have the same $\GA$-content.
\item{Secondly,} this common $\GA$-content is equal to the Kronecker invariant
$t^0_n$ of $K\GA$:
$$|D_n|_{\GA}=|D_n'|_{\GA}=t^0_n(K\GA)$$
\end{description}
Applying this to the particular case $\GA_4$, we obtain
$$|D_3|_{\GA_4}=2+1-1-1=1=|D_3'|_{\GA_4},\mbox{ so }t^0_3(K\GA)=1$$
while all other $t^0_n(K\GA_4)$ vanish.

This exhausts the information provided by combining Part One of Th.\ref{th:main}
with our knowledge of $\mathcal{D}(\GA_4)$. We next turn to Part Two of
Th.\ref{th:main}; applying this  to the result obtained above, that the
contraction-stabilization of $\GA_4$ is the graph consisting of a single loop
on one vertex, we see that all $t_n(\GA_4)$ are 0, and that the only indecomposable
summand of $K\GA_4$ of type $S(p^n)$, is $S(X-1)$ with multiplicity 1.

This leaves only one further indecomposable type to examine in
relation to $K\GA_4$, namely the exceptional type $^0T^0_0$. To
sum up our work on $\GA_4$ up to this point, we have obtained for
$K\GA_4$ the following non-zero multiplicities:
$$^0t^0_2=3,\;^0t^0_3=2,\;^0t_2=1,\;t^0_3=1,;[K\GA_4:S(X-1)]=1$$
with all other multiplicities (except possibly that of $^0T^0_0$)
being 0. Note that these computations indeed check with
eqn.(\ref{check}):
$$\mathcal{V}(\GA_4)=18=2\mathcal{V}(^0t^0_2)+3\mathcal{V}(^0t^0_3)+\mathcal{V}(^0t_2)+
\mathcal{V}(t^0_3)+ 1\cdot(\mbox{deg}(X-1))=2(3)+3(2)+2+3+1$$

Applying Part Three of Th.\ref{th:main}, in the form of eqn.(20),
and utilizing the previously obtained and checked results for the
other Kronecker multiplicities for $K\GA_4$, together with the
observation that $\GA_4$ has 23 edges, we obtain
$$^0t^0_0(\GA_4)=23-1\cdot(2+3)-0-(3)(3)-(4)(2)-(1)(1)=0$$

\medskip
Thus, finally, we obtain the Kronecker decomposition
$$K\GA_4\cong 3(^0T^0_2)\oplus 2(^0T^0_3)\oplus
T_3^0\oplus ^0\!T_2\oplus S(X-1)\;.$$

\medskip
\noindent The reader is invited to check similarly, using Th.5.4,
 that the Kronecker decompositions for the binary relations
explained in Examples 2.1, 2.2 and 4.4 are:
$$K\GA_1\cong ^0\!\!T^0_1\oplus T^0_1\oplus S(X-1),\;\;
K\GA_2\cong  ^0\!\!T^0_0\oplus ^0\!\!T^0_1\oplus S(X-1),
$$
and
$$K\GA_3\cong ^0\!\!T_1\oplus T^0_1\oplus S(X^2-1)$$
(where the first of these decompositions agrees with the result
(\ref{decomp:G1}) obtained earlier by other methods.)

\end{document}